\documentclass[psamsfonts]{amsart}
\usepackage{mymacros}

\title{Notes on Factorization Algebras and TQFTs}
\date{}
\author{Araminta Amabel}

\begin{document}

\begin{abstract}
These are notes from talks given at a spring school on 
topological quantum field theory in Nova Scotia during May of 2023.
The aim is to introduce the reader to the role of factorization algebras
and related concepts in field theory. 
In particular, we discuss the relationship between
factorization algebras, $\bb{E}_n$-algebras, vertex algebras,
and the functorial perspective on field theories.
\end{abstract}
 
\maketitle

\tableofcontents

\section{Introduction}

Ideas motivated by quantum field theory connect a wide array of 
mathematics. 
From probability theory to knot theory, or differential geometry to infinity categories, 
modeling the physical aspects of field theories and their quantization
provides a wealth of interesting concepts.
Because they permeate an expanse of topics, 
it can be hard to move between perspectives without reading a large amount of background material.
The intention behind these notes is to provide an accessible
and introductory account of several different,
but mathematically close, 
descriptions of field theories.

The viewpoints chosen are those most accessible to an algebraic topologist.
They are:

\begin{itemize}
\item factorization algebras of observables,
\item $\bb{E}_n$-algebras,
\item the functorial perspective, and
\item vertex algebras.
\end{itemize}
Roughly, these are all methods for accounting for the different measurements one can make on a physical system. 

We will give the most attention to factorization algebras.
Factorization algebras are algebraic gadgets whose 
multiplicative structures are controlled by disjoint open subsets of a manifold;
see Definition \ref{def-FA}.
Historically, similar structures were first considered
under the name of chiral algebras by Beilinson and Drinfeld \cite{BD}.

Whereas factorization algebras highlight the geometry 
of the manifold governing their multiplication,
$\bb{E}_n$-algebras have a more topological flavor.
Vertex algebras are more common in the representation theory literature
and the functorial perspective is useful when looking for structural results.  
Each perspective has its benefits which we try to highlight throughout these notes.
Lastly, we end these notes with a discussion of more recent advances and open areas of research.

As these lectures were given over the course of five days,
they are hopefully a quick and informal read. 
Exercises with hints are included.

\subsection{Assumed Background}
The reader is assumed to be familiar with the language of (ordinary, not infinity) categories. 
No knowledge of field theory is necessary.
Roughly, the expectation is that the reader has the knowledge of an average first year graduate student with an interest in topology.

\subsection{Linear Overview}

The first part includes a preliminary definition
of a classical field theory
and introduces the notion of a factorization algebra.
The main takeaway of this part 
is that field theories can be studied by 
considering their algebras of observables.
We end the first part with a discussion of quantization 
in terms of observables.
In the second part, 
we focus on topological field theories
and their relationship to $\bb{E}_n$-algebras.
The third part is on 
the global invariant of observables called factorization homology.
Afterwards, in the fourth part, we discuss the Atiyah-Segal functorial 
perspective of field theories using bordism categories.
In the last part, we consider holomorphic field theories and 
their relationship to vertex algebras.
This last part also contains remarks 
on the cobordism hypothesis, 
Koszul duality,
and the Stolz-Teichner program.
After each part, exercises are included.
The exercises vary in difficulty level, depending on the reader's background.
At the end of these notes, 
hints or solutions to each exercise are given.

\subsection{Further Reading}
These notes are short, and not all-inclusive. 
Many important topics have been left out. 
The lectures these notes are based on where given in 
conjunction with two other lecture series:
one on spectral sequences by Debray 
and one on fusion categories by Delaney. 
Notes for these talks can be found online;
see \cite{Debray, Delaney}.

For more content and and in-depth discussion of the 
ideas in these notes, 
we make the following recommendations:

\begin{itemize}
\item For factorization algebras in field theory: the original reference, 
which is wonderfully readable, is the two part series \cite{CG1} and \cite{CG2}.

\item For $\bb{E}_n$-algebras: \cite{May} is a canonical reference. In the $\infty$-setting, see \cite{HA}.

\item For factorization homology: the various works of Ayala-Francis including \cite{AFtop} and \cite[Ch. 2]{Handbook} 
as well as Ginot's notes \cite{Ginot}.

\item For the cobordism hypothesis: see \cite{FreedCob} for an overview, \cite{LurieCobordism} for Lurie's proof sketch, \cite{SP} for the details of the 2-dimensional case, and \cite{Sche} for a construction of a target category.

\item For vertex algebras: \cite{BenZviFrenkel} or \cite{Kac}.

\end{itemize}

\subsection{Acknowledgements}
As these notes were written for the Atlantic TQFT Spring School in May of 2023,
the author would like to thank the organizers 
Theo Johnson-Freyd and Geoff Vooys.
Thanks are also due to the participants of the Atlantic TQFT Spring School
for their questions. 
In particular, the author thanks Will Stewart, Eilind Karlsson, and Mathew Yu 
for serving as the TAs for the spring school.
We thank Owen Gwilliam and Hiro Lee Tanaka for providing comments and corrections 
on an earlier draft. 
The author was
supported by NSF grant DMS-1902669.

\newpage
\part{Day 1}

I am supposed to teach you about factorization algebras. 
There is a lot to say about these algebraic gadgets,
so I will not come close to covering everything.
I hope you come away with three skills:

\begin{enumerate}
\item knowing what a factorization algebra \emph{is} (the definition and examples),
\item being able to relate them to other notions of field theories, and
\item a curiosity for learning more about the subject.
\end{enumerate}

\noindent The five lectures will roughly cover the following themes:

\begin{enumerate}
\item[\textit{Day 1}.] physical motivation
\item[\textit{Day 2}.] topological field theories
\item[\textit{Day 3}.] factorization homology
\item[\textit{Day 4}.] functorial perspective
\item[\textit{Day 5}.] holomorphic field theories and applications
\end{enumerate}

For today, we have the following goal.

\begin{goal}
Discover ``factorization algebras" in nature.
\end{goal}

The quotes are because I have not yet told you the definition of factorization algebras,
and part of the goal is to stumble upon the definition ourselves.
Here, nature means physics.

\begin{rmk}
The story contained here is roughly mimicking the introduction in \cite{CG1}.
We refer the reader there for more details.
\end{rmk}

\section{Classical Story}

If you have not heard of a field theory before, 
the picture I want you to have in mind is the following.
Imagine we have some container $X$ 
(a manifold or more general space)
and a particle moving around in $X$. 

For $I\subset\bb{R}$ a time interval,
the space of maps 
\[\mathrm{Map}(I,X)\]
is all the paths a particle could take.
So $t\in I$ maps to to the position of the particle at time $t$.

This system of a particle in $X$
might be subject to constraints in the real world.

\begin{ex}
Only paths of least energy.
\end{ex}

\begin{ex}
$X=\bb{R}^2$ and only straight line paths.
\end{ex}

\begin{ex}
$X=S^2$ the sphere, and only  great circle paths.
\end{ex}

\begin{rmk}
These are examples of the general case 
where paths are geodesics.
\end{rmk}

The physical constraints on what paths a particle can take 
are called the \emph{equations of motion}
 or the \emph{Euler-Lagrange equations}.
 These can be written as a PDE, 
 and are determined by a map 
 \[S\colon\mathrm{Map}(I,X)\rta\bb{R}\]
 called the \emph{action functional}.
 The allowable paths are then 
 \[\mathrm{EL}\subset \mathrm{Map}(I,X)\]
 the subset of paths so that 
 \[\mathrm{EL}=\{f\colon I\rta X:(dS)(f)=0\}.\]
 This is the \emph{critical locus} of $S$.

More generally, we might be interested in how the particle 
moves around over some time interval $I$ 
and as we move around some parameter space $N$. 
Then $N\times I=M$ is ``spacetime" 
and we are interested in parameterized paths
\[\mathsf{Map}(M,X).\]

For now we will take the following definition.

\begin{defn}
A \emph{classical field theory} is a space of fields 
$\mathsf{Map}(M,X)$ 
and a set of equations of motion that can be encoded by the critical locus of a map
\[S\colon\mathsf{Map}(M,X)\rta\bb{R}.\]
The \emph{dimension} of the field theory is the dimension of $M$.
\end{defn}

\begin{rmk}
The action functional $S$ must be local: it must arise as the integral over $M$
of some polynomial in a field $\phi$ and its derivatives;
see \cite[Ch. 2. \S7]{CG1}.
\end{rmk}

\begin{rmk}
We already have several equivalent descriptions of a given classical field theory.
One could hand you the space of fields 
and the action functional,
or the equations of motion.
For more equivalent descriptions, 
see \cite[Ch. 4]{CG2}
and in particular \cite[Def. 4.2.0.4]{CG2} and \cite[Def. 4.4.0.5]{CG2}.
\end{rmk}

\begin{quest}
What if someone just told you the critical locus $\mathrm{EL}$?
Would that be equivalent data to the whole field theory?
\end{quest}

For an answer, see \cite[Def. 4.2.0.4]{CG2}.

\begin{ex}[Classical mechanics, massless free theory]
Say $I=[a,b]$. 
Take fields 
\[\mathsf{Map}(I,\bb{R}^n).\]
The action functional sends a field $f$ to
\[S(f)=\int_a^b\langle f(t),\frac{d^2}{dt^2}f(t)\rangle dt.\]
Here, we are taking the inner product on $\bb{R}^n$. 
In this case,
the critical locus is straight lines,
\[\mathrm{EL}=\{\text{straight lines}\}.\]
For details see \cite[Ch. 1 \S1]{CG1}.
If we add in a nonzero mass,
the free theory with mass $m$ is described in \cite[Ch. 2 \S2]{CG1}.
\end{ex}

\begin{quest}
Can you extend this notion to a massless free theory 
on a general $n$-manifold $M$?
What extra structure will you need $M$ to have?
\end{quest}

For an answer, see or \cite[\S 4.5.1]{CG2}.

\begin{ex}[Gauge Theory]
Given a Lie group $G$ and a spacetime $M$,
$G$ gauge theory on $M$ 
which I will denote $\mathsf{Gauge}_G^M$ 
has fields 
\[\mathsf{Map}(M,B_\nabla G)=\mathsf{Bun}^\nabla_GM.\]
Note that principal $G$-bundles $\mathsf{Bun}_GM$ on $M$
are classified by maps into the classifying space $BG$.
Here, $B_\nabla G$ is an upgraded version of $BG$ that
classifies principal $G$-bundles with connection.
The space $B_\nabla G$ is the quotient of $\Omega^1(-;\mathfrak{g})$ 
by the adjoint action; see \cite[\S13.3]{DiffCoh}.
See \cite{FreedHopkins} for more details on $\mathsf{Bun}_G^\nabla(M)$. 
A field is then a principal $G$-bundle $P$ on $M$
with connection $A$.
There are many different choices for action functionals,
giving different gauge theories. 
\end{ex}

\begin{ex}[Yang-Mills]\label{ex-YM}
One type of gauge theory $\mathsf{Gauge}_\bb{R}^M$ 
is Abelian Yang-Mills.
The action functional is
\[S(A)=-\frac{1}{2}\int_M \mathrm{tr}(dA\wedge *dA).\]
Here $*$ is the Hodge star operator,
More generally, we can consider Yang-Mills for an arbitrary
compact Lie group $G$.
In this case, $dA$ is replaced with the curvature of $A$.
For details, see \cite[Ch. 2 \S2]{CG1}.

This is particularly studied when $G=SU(n)$ and $M$ has dimension 4.

We can alter this action functional to produce new theories by adding 
another term in the integrand.
\[S(A)=-\frac{1}{2}\int_M \mathrm{tr}(dA\wedge *dA)+\theta\int_M\mathrm{tr}(dA\wedge dA)\]
This additional term is called the ``theta term";
see \cite[Pg. 362]{WittenTQFT} or \cite[\S 2.2]{Tong}.
Note that the theta term does not depend on the connection,
and hence does not change the equations of motion. 
Because of this lack of dependence, 
the theta term is referred to as a ``topological term."
\end{ex}

\begin{quest}
Can you show that the examples of action functionals given above
are local?
\end{quest}

\begin{rmk}
The space of fields does not have to be a mapping space.
When it is, it is referred to as a \emph{sigma model}.
\end{rmk}

\section{Measurements}

Knowing a particle is moving in a box is good,
but we want to know more.

\begin{ex}
Say you hand a little kid a box containing a spider 
(i.e. a scary particle).
They look up at you really scared. 
You tell them, 

\begin{quote}
Do not worry, 
I can tell you what type of paths the spider can take.
It only moves along straight lines!
\end{quote}

They are still scared.
And they have questions.

\begin{itemize}
\item Where in the box is the spider right now?
\item How fast is the spider moving?
\item Is the box open?
\item Are there holes in the box?
\item Can the spider get out of the box?
\end{itemize}
\end{ex}

\begin{rmk}
Notice that these questions are of three different types.
The first two are about the particle and the equations of motion.
They can maybe be answered locally.
The third and fourth questions are about the topology of the spacetime.
That is global information.
Answering the last question requires knowledge about the whole field theory.
\end{rmk}

These questions are types of \emph{measurements}
one can make on the system.
In classical field theory,
one can (in theory) answer all of these questions.
This changes in quantum field theory!

\begin{thm}[Heisenberg Uncertainty Principle]
In quantum field theory,
one cannot precisely know both the position and the momentum
of a particle at the same time.
\end{thm}

\begin{goal}
Formulate this uncertainty mathematically.
\end{goal}

\subsection{Observables}
\begin{defn}
Given a classical field theory
\[\mathrm{EL}\subset\mathrm{Maps}(I,X),\]
the {classical observables} are the 
real functions on the critical locus,
\[\mathrm{Obs}^\mathrm{cl}=\mathcal{O}_{\mathrm{EL}}.\]
\end{defn}

Thus, an observable takes in an allowed path 
and spits out a real number:
the measurement on that path.

\begin{rmk}
I did not say what kind of functions we are taking.
This will depend on the context.
Generally the (derived) critical locus lives in the world of derived $C^\infty$ geometry;
see \S\ref{sec-CL}.
When things are nice, holomorphic functions also makes sense.
\end{rmk}

The vector space of functions 
\[\mathrm{Obs}^\mathrm{cl}=\mathcal{O}_\mathrm{EL}\]
form a commutative algebra;
\[f,g\colon\mathrm{EL}\rta\bb{R}\]
combine to give
\[(fg)(u)=f(u)g(u).\]
Physically,
given two measurements $f,g$,
we can perform them at the same time
to get a new measurement $fg$.

\underline{Note}: This is exactly what fails in the quantum world!

There, we \emph{cannot} take two measurements 
(like position and momentum) simultaneously. 

What can we do?
Say $I=[a,b]$ is our spacetime. 
Let $f,g$ be observables.
We can form a new measurement on $I$ 
that on 
$(t_1,t_2)$ does $f$,
on $(t_3,t_4)$ does $g$
and does nothing in between.
Nothing here means the measurement sending every path to zero.

So we have a way of combining observables 
on disjoint open intervals. 

More generally, say we were looking at a spacetime
\[M=I\times N\]
where $N$ has dimension $n$. 
We can perform different measurements on disjoint disks.
Given two embeddings 
\[\bb{D}^{n+1}\simeq I\times \bb{D}^n\rta I\times N\]
we do $f$ on one and $g$ on the other. 

\underline{Upshot}. 
Whatever quantum field theory is,
the set of measurements (or observables)
one can make on it have 
a weird multiplication structure controlled 
by disjoint open subsets of spacetime.

The resulting algebraic structure is called a \emph{factorization algebra}.

\section{Factorization Algebras}

The following can be found in  \cite[Ch. 6 Def. 1.4.2]{CG1}
and \cite[Ch. 3 \S1.1]{CG1}.

\begin{defn}\label{def-FA}
Let $M$ be a manifold.
A \emph{factorization algebra on $M$} is a functor
\[\mathcal{F}\colon\mathsf{Open}(M)\rta\mathsf{Ch}\]
together with maps
\[\mathcal{F}(U_1)\otimes\cdots\otimes\mathcal{F}(U_k)\rta\mathcal{F}(V)\]
for disjoint unions of opens $U_1,\dots, U_k$ in $M$ contained in a single open subset $V$, 
\[U_1\sqcup\cdots\sqcup U_k\subset V,\]
that are equivalences if 
\[U_1\sqcup\cdots\sqcup U_k= V,\]
and such that the maps are compatible
\[\begin{xymatrix}
{\bigotimes_{i=1}^l\bigotimes_{j=1}^{k_i}\mathcal{F}(U_{i,j})\arw[rr]\arw[dr] & & \bigotimes_{i=1}^l\mathcal{F}(V_i)\arw[dl]\\
& \mathcal{F}(W) & }
\end{xymatrix}\]
for disjoint opens $U_{i,1},\dots,U_{i,k_i}\subset V_i$ inside 
disjoint opens $V_1,\dots,V_l$ inside an open $W$.
Moreover, 
$\mathcal{F}$ satisfies a (Weiss) cosheaf condition.
\end{defn}

Let's break that all down.

A factorization algebra is the data of:

\begin{itemize}
\item for an inclusion of disjoint disks
\[\bigsqcup_{i=1}^k D_i\hookrightarrow M,\]
a vector space 
\[\mathcal{F}\left(\bigsqcup_{i=1}^k D_i\right)\]
and an isomorphism
\[\mathcal{F}\left(\bigsqcup_{i=1}^k D_i\right)\simeq\bigotimes_{i=1}^k \mathcal{F}(D_i),\]
and 
\item for an inclusion
\[\bigsqcup_{i=1}^k D_i\subset D\hookrightarrow M,\]
a map
\[\mathcal{F}\left(\bigsqcup_{i=1}^k D_i\right)\rta \mathcal{F}(D).\]
\end{itemize}

The weird multiplication comes from putting this data together:
\[\bigotimes_{i=1}^k \mathcal{F}(D_i)\simeq \mathcal{F}\left(\bigsqcup_{i=1}^k D_i\right)\rta\mathcal{F}(D)\]
tells us how to multiply observables on the disjoin disks $D_i$.

The last piece was the cosheaf condition.
This is going to be on the exercises.
A reference for Weiss covers is \cite[Ch. 6. \S 1.2]{CG1}.

Whatever a quantum field theory is,
we should be able to take measurements on it.
The set of these measurements is called the 
\emph{quantum observables}.
In keeping with the spirit of these notes,
rather than defining what a quantum field theory is,
we just describe the set of observables.

To obtain a precise mathematical description of 
a quantum field theory,
we need to assume that the QFT is perturbative.
Under this assumption we have the following.

\begin{thm}[Costello-Gwilliam]\label{thm-Cg1}
The quantum observables $\mathsf{Obs}^q$ 
of a field theory on spacetime $M$
has the structure of a factorization algebra on $M$.
\end{thm}

\noindent This is \cite[Thm. 8.6.0.1]{CG2}.

\begin{rmk}
In the nonperturbative setting,
the quantum observables form a prefactorization algebra;
see \cite[Ch. 1, \S 3]{CG1}.
\end{rmk}

\begin{rmk}
Theorem \ref{thm-Cg1}
works for QFTs in Riemannian signature. 
For a Lorentzian signature analogue,
see \cite{KasiaGwilliam}.
\end{rmk}

The factorization algebra structure on $\mathsf{Obs}^q$
means we get a vector space 
\[\mathsf{Obs}^q(U)\]
for every open subset $U$ of spacetime.
These are the \emph{local} observables at $U$.
When we were discussing classical observables,
we had a single vector space 
\[\mathsf{Obs}^\mathrm{cl}=\mathrm{Hom}(\mathrm{EL},\bb{R}).\]
One could ask how to get a single 
object $\mathsf{Obs}^q(M)$
from the factorization algebra $\mathsf{Obs}^q$.
This should be some sort of ``global" observables.
Since $\mathsf{Obs}^q$ is a type of cosheaf,
we can take its global sections.

\begin{defn}
Let $\mathcal{F}$ be a factorization algebra on $M$.
The \emph{factorization homology} of $\mathcal{F}$
is the global sections 
\[\int_M\mathcal{F}=\mathcal{F}(M).\]
\end{defn}

\section{Quantization}

Previously we talked about how classical observables form a commutative algebra,
and quantum observables form a factorization algebra.
I want to end today by talking about how
to \emph{quantize} a field theory.

To get at this question,
we need an easy, workable example.
Let $\bb{R}$ be our spacetime,
and take as fields 
\[\mathsf{Map}(\bb{R},T^*\bb{R}).\]
Lastly, 
our action functional is going to be $S=0$.
From the exercises, we know that 
\begin{itemize}
\item the critical locus $\mathrm{EL}$ of $S$ has functions
\[\mathsf{Obs}^\mathrm{cl}=C^\infty(T^*\bb{R}),\]
\item and the quantum observables are a factorization algebra on $\bb{R}$. 
Assume that $\mathsf{Obs}^q$ comes from an associative algebra $A$.
\end{itemize}

On the level of observables,
we can rephrase our quantization question as:

\begin{center}
Given a commutative algebra,
how do we get a factorization algebra out of it?
\end{center}

The process of going from a commutative algebra
to an associative algebra
is called \emph{deformation}.

\begin{defn}
Let $T$ be a commutative algebra.
A \emph{deformation} of $T$ is an associative algebra 
structure on $T[[\hbar]]$ such that 
\[T[[\hbar]]/\hbar\simeq T\]
as algebras.
\end{defn}

If this was all we asked for, 
we could always take 
\[\mathsf{Obs}^q=\mathsf{Obs}^\mathrm{cl}[[\hbar]]\]
with the usual multiplication.
That is not telling us anything about the field theory.
We need to ask for the deformation to encode more information.

\begin{ex}
In our running example, we have
\[\mathsf{Obs}^\mathrm{cl}=C^\infty(T^*\bb{R})=\bb{R}[p,q].\]
We have an interesting structure on $T^*\bb{R}$:
the symplectic form. 
On functions, the symplectic form gives a Poisson bracket,
\[\{p,q\}=1.\]
\end{ex}

\begin{rmk}
In general, the derived critical locus $\mathrm{EL}$ of $S$
has the structure of a $(-1)$-shifted symplectic stack,
so $\mathcal{O}_\mathrm{EL}$ is a $P_0$-algebra.
\end{rmk}

We can ask for deformations of classical observables 
that respect this Poisson bracket.
That is,
a deformation $R=T[[\hbar]]$ 
so that 
\[[f,g]=\hbar\{f,g\}\]
up to higher order terms in $\hbar$,
for $f,g\in T$.

We can think of this as saying that the Poisson structure is 
``semi-classical." 
It tells us the first step in the quantum direction, 
that is the $\hbar$ term.

\begin{ex}
The quantum observables of our example theory 
is the Weyl algebra
\[\mathsf{Obs}^q=\mathbb{R}[[\hbar]][p,q]\]
with multiplication so that $[p,q]=\hbar$. 
\end{ex}

More generally,
we could have a theory with fields
\[\mathrm{Map}(\bb{R},V)\]
where $V$ was a symplectic manifold.
Then $\mathcal{O}_\mathrm{EL}$ again has a Poisson bracket,
and quantum observables are a deformation respecting this bracket.

\begin{rmk}
One can in fact deform functions on any Poisson manifold,
using a globalized version of the Weyl algebra;
see \cite{Kontsevich} 
wherein one can also find a formal definition of deformation quantization.
\end{rmk}

In higher dimensions, 
we ask for the space of fields to have a shifted symplectic structure,
and use this to deform the observables.
\begin{quote}
Thus in the language of factorization algebras,
quantization is deformation.
\end{quote}

See \cite[\S 1.4]{CG2} for details.

\section{Exercises}

There are lot of exercises here of varying difficulty. 
You do not need to do all of them,
or any of them completely. 
Just get a feel for things and have fun learning!

\subsubsection{Cosheaf Condition}

A good reference for these questions is 
\cite[Ch. 6 \S 1]{CG1}.

Recall that a functor
\[F\colon\mathsf{Open}(M)\rta\mathcal{C}\]
is a \emph{cosheaf}
if for every open cover $\mathcal{U}$ of $M$,
the diagram
\[\bigsqcup_{i,j}F(U_i\cap U_j)\rightrightarrows \bigsqcup_k F(U_k)\rta F(M)\]
is a colimit diagram.

\begin{defn}
An open cover $\{U_j\}$ of $M$ 
is a \emph{Weiss} cover if for every finite set $x_1,\dots,x_k$ 
of points in $M$,
there exists some $U_j$ so that 
\[\{x_1,\dots, x_k\}\subset U_j.\]
\end{defn}

For a detailed introduction to cosheaves,
see \cite[Part 1]{Curry}.

\begin{quest}\label{q-5.new0}
A canonical example of a sheaf is smooth functions on a manifold $M$. 
More generally, given a vector bundle $\pi\colon V\rta M$,
show that sections of $\pi$ form a sheaf.
In contrast, show that compactly supported sections of a vector bundle form a cosheaf.
\end{quest}

\begin{quest}\label{q-5.2}
Given an example of an open cover (in the usual topology) of $\bb{R}^2$ 
that is not a Weiss cover.
\end{quest}

\begin{quest}\label{q-5.3}
Give an example of a Weiss cover that exists on any manifold $M$.
\end{quest}

\begin{defn}
Let $\mathsf{Ran}(M)$ be the set of nonempty finite subsets of $M$.
\end{defn}

Take on faith for a moment that the set $\mathsf{Ran}(M)$ has a reasonable topology.
For details, see \cite[\S5.5.1]{HA}
and \cite{Cepek}.

\begin{quest}\label{q-5.5}
Show that a Weiss cover for $M$ determines an ordinary cover for $\mathsf{Ran}(M)$.
Conversely,
show that an ordinary cover for $\mathsf{Ran}(M)$ determines an Weiss cover on $M$. 
\end{quest}

\noindent Ignore topology issues for this, 
just show it on the level of sets.

It turns out that $\mathsf{Ran}(M)$ has a topology
so that cosheaves on $\mathsf{Ran}(M)$ are the same as
cosheaves on $M$ for the Weiss topology;
see \cite[Rmk. 25]{Ginot} and the references therein.

\subsubsection{Examples}

\begin{quest}\label{q-1comm}
Given a field theory with space time $M$, 
classical observables are a commutative algebra. 
Show that classical observables form a factorization algebra on $M$.
In fact, any commutative algebra forms a factorization algebra on $M$.
\end{quest}

\begin{quest}\label{q-5.7}
Show that an associative algebra forms a factorization algebra.
\end{quest}

Note that I did not tell you what space on which it should be a factorization algebra;
this is part of the exercise.
Call this space $X$.

\begin{quest}\label{q-5.8}
Let $X$ be the space on which associative algebras determine a factorization algebra.
When is a factorization algebra $\mathcal{F}$ on $X$
determined by an associative algebra?
\end{quest}

\begin{quest}\label{q-5.9}
Let $\mathcal{F}$ be a factorization algebra on $[0,1]$
so that 
\[\mathcal{F}\vert_{(0,1)}\]
comes from an associative algebra. 
Describe the structure of $\mathcal{F}$ over $[0,1]$ 
in more familiar terms. 
\end{quest}

\subsubsection{Critical Locus}\label{sec-CL}

Recall that the classical observables are 
functions on the critical locus,
\[\mathsf{Obs}^\mathrm{cl}=\mathcal{O}_{\mathrm{EL}}\]
where 
\[\mathrm{EL}\subset\mathsf{Map}(M,X)\]
is that the critical locus of the action functional
\[S\colon\mathsf{Map}(M,X)\rta\bb{R}.\]

\begin{quest}\label{q-5.new}
Take $M$ to be a point and $X=W$ to be a vector space.
Describe the critical locus of an action functional
\[S\colon\mathsf{Map}(\mathrm{pt},W)\rta\bb{R}.\]
\end{quest}

For convenience, 
let us now pretend the mapping space $\mathrm{Map}(M,X)$
(the fields) is a smooth manifold $Y$.

\begin{quest}\label{q-5.10}
Let $\Gamma(dS)\subset T^*Y$
denote the graph of $dS$. 
Show that the critical locus of $S$ 
is the intersection 
\[\Gamma(dS)\cap \mathrm{Zero}_Y\]
where $\mathrm{Zero}_M$ is the zero-section of $Y$ in $T^*Y$.
\end{quest}

\begin{quest}\label{q-5.11}
What is $\mathrm{EL}$ if $S=0$?
\end{quest}

In reality, we want $\mathrm{EL}$ to be a fancier version of the critical locus:
the \emph{derived} critical locus. 
The derived critical locus of $S\colon Y\rta \bb{R}$ is a dg space on $Y$,
so it is determined by its functions (which are a chain complex).

\begin{defn}
The derived critical locus of $S\colon Y\rta \bb{R}$ has functions the derived tensor product
\[C^\infty(\Gamma(dS))\bigotimes_{C^\infty(T^*Y)}^\bb{L}C^\infty(Y).\]
\end{defn}

See \cite[\S 4.1 and \S4.3]{CG2} for details.

\begin{quest}\label{q-5.13}
What is the derived critical locus if the space of fields is 
\[\mathsf{Map}(\mathrm{pt},W)\]
as in Question \ref{q-5.new}?
What is the derived critical locus of $S=0$?
\end{quest}

\begin{quest}\label{q-5.14}
What is the derived critical locus of general $S$ in terms of $S$?
\end{quest}

\newpage
\part{Day 2}

\section{TQFTs}

Yesterday we talked about how classical observables form a commutative algebra,
and quantum observables form a factorization algebra.
Today I want to talk about how this structure behaves when
we have a \emph{topological} field theory.

Informally, a field theory is topological if 
it does not depend on a metric.
Let's investigate what this means for quantum observables.
Let $M$ be our spacetime.
If our theory is topological the cosheaf
\[\mathsf{Obs}^q\colon\mathsf{Open}(M)\rta\mathcal{C}\]
cannot know about measurements of size or distance.
For example, consider the value
$\mathsf{Obs}^q(B_r(0))$
on a ball of radius $r$. centered at the origin.
Since $\mathsf{Obs}^q$ does not know about size,
it cannot distinguish between balls of different radii. 
Thus 
\[\mathsf{Obs}^q(B_r(0))\xrta{\sim}\mathsf{Obs}^q(B_{r'}(0))\]
for $r<r'$. 
Moreover, it does not know distance from the origin.
So if we move $B_r(0)$ around inside $B_{r'}(0)$,
that will not effect the answer either.

\begin{defn}
A factorization algebra $\mathcal{F}$ on $M$ is \emph{locally constant}
if for any inclusion of disks 
\[D_1\subset D_2\]
in $M$,
the induced map
\[\mathcal{F}(D_1)\rta\mathcal{F}(D_2)\]
is an equivalence.
\end{defn}
This is \cite[Ch. 3 Def. 2.0.1]{CG1}.

Motivated by . \cite[Ch. 1 \S 4.3]{CG1}, 
we are going to use this as a definition.

\begin{defn}
A field theory is \emph{topological} if $\mathsf{Obs}^q$ is locally constant.
\end{defn}

\begin{ex}
A locally constant factorization algebra on $\bb{R}$ is an associative algebra.
\end{ex}

\begin{ex}
A locally constant factorization algebra on $\bb{R}^\infty$ is a commutative algebra.
\end{ex}

\begin{ex}
The data of a locally constant factorization algebra on $\bb{R}^2$ is 
a vector space $V$ with many types of multiplications and coherencies. 
\end{ex}

Locally constant factorization algebras on euclidean spaces 
give a family of algebra structures starting from associative,
and becoming more commutative as the dimension increases.

To precisely describe this structure, 
we will use the language of operads.

\section{Operads}

Informally, 
an \emph{operad} $\mathcal{Q}$ is a sequence of spaces $\mathcal{Q}(k)$ encoding $k$-ary operators of an algebraic structure. 
For example, there is an operad $\msf{Assoc}$ that records the data of a unit, a multiplication map, and the associativity conditions for an associative algebra. 
To formalize this, we need some preliminary definitions.

Let $\msf{Fin}^\mrm{bij}$ be the category of finite sets and bijections. 
The category of \emph{symmetric sequences in} $\msf{Spaces}$ 
is the functor category 

\[\msf{Sseq}(\msf{Spaces})\colonequals \msf{Fun}(\msf{Fin}^\mrm{bij},\msf{Spaces}).\]

Symmetric sequences can be given the structure of a monoidal category as follows. 
The composition product $R\circ S$ of two symmetric sequences is

\begin{align*}
(R\circ S)(n)=\bigoplus_i R(i)\otimes_{\Sigma_i}\left(\bigotimes_{j_1+\cdots+j_i=n} (S(j_1)\otimes\cdots\otimes S(j_i))\times_{\Sigma_{j_1}\times\cdots\times\Sigma_{j_i}}\Sigma_n\right). 
\end{align*}

The unit of the composition product, denoted $\mathcal{O}_{\msf{triv}}$, sends a finite set $B$ to the unit $\bb{1}_\msf{Spaces}$ of $\msf{Spaces}$ if $|B|=1$ and to the zero object $*$ of $\msf{Spaces}$ otherwise. 

\begin{defn}
An \emph{operad} in $\msf{Spaces}$ is a monoid object in symmetric
sequences  $\msf{Sseq}(\msf{Spaces})$.
\end{defn} 

An operad $\mathcal{O}$ in spaces has an underlying functor $\msf{Fin}^\mrm{bij}\rta\msf{Spaces}$. 
For each $k\in\bb{N}$, we denote by $\mathcal{O}(k)$ the image of the finite set with $k$ elements $[k]$ under this functor.
This matches with \cite[Def. 1.1]{May}.

\begin{ex}[Little $n$-Disks Operad]
Define an operad $\bb{E}_n$ in spaces by 
\[\bb{E}_n(k)=\mathrm{Conf}_k(\bb{R}^n),\]
the configuration space of $k$ distinct points in $\bb{R}^n$, 
topologized as a subset of $(\bb{R}^n)^k$.
The easiest way to see the product 
\[\bb{E}_n\circ\bb{E}_n\rta\bb{E}_n\]
is to consider each point in $\bb{R}^n$ as a little open disk centered at that point.
Then the product is just the inclusion of disks.
For more details, see \cite[Ch. 4]{May}.
A good reference for configuration spaces in general is \cite{KnudsenConf}.
\end{ex}

We think of the space $\mathcal{O}(k)$ 
as parameterizing $k$-ary operations for some type of algebraic structure.

\begin{defn}
An $\mathcal{O}$-algebra in $\mathcal{C}$ 
is an object $V\in\mathcal{C}$ together with maps 
\[\mathcal{O}(k)\otimes V^{\otimes k}\rta V\]
for every $k$, 
compatible with the multiplication maps for $\mathcal{O}$.
\end{defn}

Here, the tensor product $\mathcal{O}(k)\otimes V$ 
is tensoring together a space $\mathcal{O}(k)$
and an object $V\in \mathcal{C}$. 
To make sense of this in general, 
$\mathcal{C}$ must be tensored over spaces.
This is true if, for example, $\mathcal{C}$ is 
an $\infty$-category; see \ref{S-infty}. 
For now, we have the following examples.

\begin{itemize}
\item For $\mathcal{C}=\mathsf{Sets}$, we tensor $V$ with the set $\pi_0(\mathcal{O}(k))$.

\item For $\mathcal{C}=\mathrm{Vector}$, we tensor $V$ with the vector space $H_0(\mathcal{O}(k))$.

\item For $\mathcal{C}=\mathsf{Ch}$, we tensor $V$ with the chain complex $C_\bullet(\mathcal{O}(k))$.

\item For $\mathcal{C}=\mathsf{Cat}$, we tensor $V$ with the fundamental groupoid $\pi_\bullet(\mathcal{O}(k))$.
\end{itemize}

\begin{thm}[Lurie]
Locally constant factorization algebras on $\bb{R}^n$ are the same as 
$\bb{E}_n$-algebras.
\end{thm}

\noindent This is \cite[Thm. 5.5.4.10]{HA}.

\subsection{Factorization Homology}

Recall that for a factorization algebra $\mathcal{F}$ 
on $M$,
we had a notion of \emph{factorization homology}
given by global sections 
\[\int_M\mathcal{F}=\mathcal{F}(M).\]

We want to investigate this local to global structure in the locally constant case.
For this,
note that there is an equivalence 
\[\mathsf{Emb}^\mathrm{fr}\left(\bigsqcup_k\bb{R}^n,\bb{R}^n\right)=\mathsf{Conf}_k(\bb{R}^n)\]
which records the center of the disks.

We are going to construct a notion of factorization homology
of an $\bb{E}_n$-algebra $A$ over a framed $n$-manifold $M$;
\[\int_MA.\]

\section{Digression on Homotopy Coherence}\label{S-infty}

The target category for our factorization (or $\bb{E}_n$) algebras
has been chain complexes $\mathsf{Ch}$.
As we will explain here, 
to distinguish between $\bb{E}_n$-algebras and commutative algebras
for $n>1$, 
we need $\mathsf{Ch}$ to be the \emph{$\infty$-category} of chain complexes.
More generally, 
there are interesting notions of factorization (or $\bb{E}_n$) algebras valued 
in any symmetric monoidal $\infty$-category $\mathcal{C}$. 

\begin{rmk}
Later in these notes, everything will start to become $\infty$-categories
for things to be true as stated.
If you do not care for category theory,
ignore the rest of this section.
\end{rmk}

For a reference on $\infty$-categories see \cite{HTT} for detail
or \cite{LurieInf} for an introduction.
See also \cite[Ch. 2]{Tanaka} for an introduction to $\infty$-categories
in the factorization homology context.

If you are not familiar with $\infty$-categories,
you can get through these notes by keeping three ideas in mind:

\begin{itemize}

\item Chain complexes up to quasi-isomorphism $\mathsf{Ch}$ 
is a good example of an $\infty$-category.

\item In an $\infty$-category, there are \emph{spaces} of morphisms rather than sets of morphisms.

\item To actually say things, like writing a functor between $\infty$-categories, 
extra care must be taken. Rather than writing down what a map does on objects and morphisms, things like universal properties must be used.

\end{itemize}

Our main example of an $\infty$-category will be $\mathsf{Ch}$. 

Ordinary categories are special types of $\infty$-categories.
For example, the category of real vector spaces $\mathrm{Vector}$
or the category of chain complexes $\mathrm{Chain}$. 

One can also build more interesting $\infty$-categories from
ordinary categories.
Let $D$ be an ordinary abelian category.
Form the (ordinary) category $\mathrm{Chain}(D)$ of chain complexes in $D$.
Regard $\mathrm{Chain}(D)$ as an $\infty$-category.
The \emph{derived $\infty$-category} of $D$ 
is obtained from $\mathrm{Chain}(D)$ by inverting quasi-isomorphisms in the land 
of $\infty$-categories. 
For a reference on derived categories see \cite{Keller}.

\begin{rmk}
If one inverts the quasi-isomorphisms of $\mathrm{Chain}(D)$ regarded as an
ordinary category and \emph{not} in the land of $\infty$-categories,
one obtains the homotopy category of the derived $\infty$-category of $D$.
\end{rmk}

For example, the derived $\infty$-category of $\mathrm{Vector}$ is 
$\mathsf{Ch}$. 
The notation $\mathsf{Vect}$ is sometimes used to denoted the derived $\infty$-category of vector spaces. 
We will stick to the notation $\mathsf{Ch}$ here. 

Similarly, the homotopy theory of spaces
by inverting weak equivalences in the category of topological spaces.
The resulting $\infty$-category is denoted $\mathsf{Spaces}$ in these notes.
We also get the $\infty$-category of ordinary categories (denoted $\mathsf{Cat}$)
by inverting equivalences of categories. 

\begin{prop}
Let $\mathcal{C}$ be an ordinary symmetric monoidal category, such as $\mathrm{Vect}$. 
Then the forgetful functor
\[\mathsf{Alg}_{\bb{E}_\infty}(\mathcal{C})\rta\mathsf{Alg}_{\bb{E}_n}(\mathcal{C})\]
is an equivalence for all $n>1$.
\end{prop}

\begin{proof}
This follows from the Eckmann-Hilton theorem \cite{Eckmann}
which roughly says that a set with two unital binary operations 
which commute is a commutative algebra.
If $A$ is a $\bb{E}_n$-algebra for $n>1$,
then there are multiple compatible binary operations 
from the different configurations $\mathsf{Conf}_2(\bb{R}^n)$.
\end{proof}

\section{Exercises}

\begin{quest}\label{q-10.1}
Show that a locally constant factorization algebra on $\bb{R}^n$
determines a locally constant factorization algebra on $\bb{R}^{m}$
for any $m<n$.
\end{quest}

\begin{quest}\label{q-10.2}
Write out and try to visualize the multiplication  
\[(\bb{E}_n\circ\bb{E}_n)(k)\rta\bb{E}_n(k)\]
for some small values of $n$ and $k$
such as $n=1,2$ and $k=1,2,3$. 
\end{quest}

\begin{quest}\label{q-10.3}
Show that $\bb{E}_1$-algebras in vector spaces are associative algebras.
\end{quest}

\begin{quest}\label{q-10.4}
Show that $\bb{E}_\infty$-algebras in vector spaces are commutative algebras.
\end{quest}

\begin{quest}\label{q-10.6}
Show that an $\bb{E}_1$-algebra in $\mathsf{Cat}$
is a monoidal category.
\end{quest}

\begin{quest}\label{q-10.7}
Show that an $\bb{E}_2$-algebra in $\mathsf{Cat}$
is a braided monoidal category.
\end{quest}

\begin{quest}\label{q-10.5}
Show that an $\bb{E}_\infty$-algebra in $\mathsf{Cat}$
is a symmetric monoidal category.
In fact, $\bb{E}_n$-algebras in $\mathsf{Cat}$ for $n>2$ 
are symmetric monoidal categories.
\end{quest}

\begin{quest}\label{q-10.8}
Let $X$ be a pointed topological space.
Show that $\Omega^nX$ is an $\bb{E}_n$-algebra in spaces.
\end{quest}

In fact, every (connected, grouplike) $\bb{E}_n$-algebra in spaces
looks like $\Omega^nX$ for some $X$.
This is called May's Recognition Principle, \cite[Thm. 1.3]{May}.

\subsubsection{Enveloping Algebras}
A good reference for these exercises is \cite[Ch. 3 \S 4]{CG1}.

For the following exercise,
recall the Chevalley-Eilenberg complex of a Lie algebra
\[C_\bullet(\mathfrak{h})=(\mathrm{Sym}(\mathfrak{h}[1]),d)\]
where $d$ is determined by the bracket on $\mathfrak{h}$;
see, for example, \cite[Appendix A Def. 3.1.2]{CG1}.

Let $\mathfrak{g}$ be a Lie algebra over $\bb{R}$.
Let $\mathfrak{g}^\bb{R}$
be the cosheaf valued in chain complexes on $\bb{R}$ assigning
\[\Omega^*_c(U)\otimes\mathfrak{g}\]
with differential $d_{dR}$ 
to an open interval $U\subset \bb{R}$. 

\begin{quest}\label{q-10.9}
Show that the assignment
\[U\mapsto H^\bullet(C_\bullet(\mathfrak{g}^\bb{R}(U)))\]
defines a factorization algebra on $\bb{R}$.
Call it $\mathcal{U}(\mathfrak{g})$.
\end{quest}

\begin{quest}\label{q-10.10}
Show that $\mathcal{U}(\mathfrak{g})$ is locally constant.
\end{quest}

\begin{quest}\label{q-10.11}
Show that $\mathcal{U}(\mathfrak{g})$ is determined by the 
associative algebra $U\mathfrak{g}$,
the universal enveloping algebra.
\end{quest}

\newpage
\part{Day 3}

\section{Factorization Homology}

Our next goal is to define factorization homology in the language 
of $\bb{E}_n$-algebras.
To do this, we will follow \cite{AFtop}.
Another reference is \cite[Def. 5.5.2.6]{HA}.

\begin{defn}
Let $\cal{M}\msf{fld}_n$ be the $\infty$-category of $n$-manifolds. 
This has objects $n$-manifolds and morphisms smooth embeddings.

Let $\cal{D}\msf{isk}_n\subset\cal{M}\msf{fld}_n$ be the full $\infty$-subcategory consisting of manifolds isomorphic to finite disjoint unions of Euclidean spaces.
\end{defn}

See \cite[Def. 2.1]{AFtop} for more details on these definitions, 
including a necessary technical condition on the manifolds:
that they admit finite good covers.

\begin{rmk}
If you are pretending these are just topological categories, the mapping spaces $\msf{Emb}(M,N)$ are given the compact-open topology.
\end{rmk}

\begin{rmk}
We consider the empty set $\empty$ to be an object of $\cal{M}\msf{fld}_n$ for every $n$. 
\end{rmk}

\begin{rmk}
We will need the following variation of disk categories.
\begin{itemize}
\item Framed disks: $\cal{D}\msf{isk}^\mrm{fr}_n$ 
will have the same objects of $\cal{D}\msf{isk}_n$ 
but with framed embeddings as morphisms. 
A framed embedding is an embedding $M\rta N$ 
so that the given framing on $M$ and 
the pulled back framing commute up to a chosen homotopy. 
\end{itemize}
\end{rmk}

Note that $\cal{D}\msf{isk}_n$ 
is a symmetric monoidal $\infty$-category 
with tensor product given by disjoint union. 
Given an $n$-manifold $M$, let $\cal{D}\msf{isk}_{n/M}$ 
denote the over category. 
Objects of $\cal{D}\msf{isk}_{n/M}$ 
are embeddings $U\hookrightarrow M$ 
so that $U$ is isomorphic to a finite disjoint union 
$\sqcup\bb{R}^n$ of Euclidean spaces. 
Morphisms in this category are triangles

\[\xymatrix{
U\arw[r]\arw[d] & M\\
V\arw[ur]	
}\]
that commute up to a chosen isotopy. 
The over category $\cal{D}\msf{isk}_{n/M}$ 
comes with a forgetful functor 
\[\cal{D}\msf{isk}_{n/M}\rta\cal{D}\msf{isk}_n\]

The following is \cite[Def. 3.1]{AFtop}.

\begin{defn}
An $n$-disk algebra $A$ with values in 
a symmetric monoidal $\infty$-category $\cal{V}$ 
is a symmetric monoidal functor 

\[A:\cal{D}\msf{isk}_n\rta\cal{V}\]
Let $\msf{Alg}_n(\cal{V})$ denote the category of $n$-disk algebras. 
\end{defn}

\subsection{$\cal{E}_n$-algebras}

If we redo everything above with framed manifolds, framed embeddings, and such, then a 
\emph{framed} $n$-disk algebra is the same as an $\cal{E}_n$-algebra. 
The equivalence goes as follows. 
Given a framed $n$-disk algebra $A$ with values in $\cal{V}$, 
define an $\cal{E}_n$-algebra in $\cal{V}$ by $A(\bb{R}^n)$ 
and action
\[\msf{Emb}\left(\coprod_I\bb{R}^n,\bb{R}^n\right)\otimes A(\bb{R}^n)^{\otimes I}\rta A(\bb{R}^n)\]
by identifying 
$A(\bb{R}^n)^{\otimes I}\simeq A(\coprod\limits_I\bb{R}^n)$ 
and applying the given embedding.

More precisely, there is an equivalence of categories
\[\msf{Alg}_{\cal{D}\msf{isk}_n^\mrm{fr}}(\cal{V})\cong\msf{Alg}_{\cal{E}_n}(\cal{V})\]
For details, see \cite[Rmk. 2.10]{AFtop}.

\begin{defn}
Let $M$ be an $n$-manifold and $A$ an $n$-disk algebra 
valued in $\cal{V}$. 
The \emph{factorization homology of $M$ with coefficients in $A$} 
is the homotopy colimit
\[\msf{colim}\left(\cal{D}\msf{isk}_{n/M}\rta\cal{D}\msf{isk}_n\xrta{A}\cal{V}\right).\]
\end{defn}

This is \cite[Def. 3.2]{AFtop}.

\begin{rmk}
The symbol $\int$ is also used in category theory to denote (co)ends.
Ordinarily, a coend over a category $\mathcal{C}$ is denoted by $\int_\mathcal{C}$
and an end is denoted $\int^\mathcal{C}$.
Although the variable $M$ appears in the bottom of the $\int$ symbol, factorization homology can be expressed as a coend
-over the symmetric monoidal envelope of $\bb{E}_n$.
See \cite[Rmk. 3.3.4]{Zero} for details.
\end{rmk}

\begin{rmk}
Factorization homology for framed manifolds can also be 
described as a bar construction in the language of operads;
see \cite[Appendix A]{My1}.
\end{rmk}

\subsection{Homology Theories for Manifolds}

Factorization homology satisfy a version, 
more suited to manifolds, 
of the Eilenberg-Steenrod axioms for homology theories. 
The main axiom of such theories is called ``$\otimes$-excision."
The following is \cite[Def. 3.15]{AFtop}.

\begin{defn}
	A symmetric monoidal functor 
	\[F:\msf{Mfld}_n\rta\cal{C}\msf{h}_\bb{Q}\]
	satisfies $\otimes$-excision if, for every collar-gluing $U\bigcup_{V\times\bb{R}}U'\simeq W$, the canonical morphism
	\begin{align}\label{eq-11}
	F(U)\bigotimes_{F(V\times\bb{R})}F(U')\rta F(W)
	\end{align}
	is an equivalence. 
\end{defn}

One reason we have restricted to collar-gluings 
is so that this tensor product makes sense.
As an exercise for this section, 
you will check that the pieces of (\ref{eq-11})
have the correct algebraic structure to form the tensor product.

\begin{defn}
	The $\infty$-category of \emph{homology theories for $n$-manifolds} 
	valued in $\cal{C}\msf{h}_\bb{Q}$ is the full $\infty$-subcategory
	\[\bb{H}(\cal{M}\msf{fld}_n,\cal{C}\msf{h}_\bb{Q})\subset\msf{Fun}^{\otimes}(\cal{M}\msf{fld}_n,\cal{C}\msf{h}_\bb{Q})\]
	of symmetric monoidal functors that satisfy $\otimes$-excision.
\end{defn}

Not only is factorization homology a homology theory for $n$-manifolds, 
it also is the only such thing.
The following is \cite[Thm. 3.26]{AFtop}.

\begin{thm}[Ayala-Francis]\label{ES}
	There is an equivalence
	\[\int:\msf{Alg}_n(\cal{C}\msf{h}_\bb{Q})\leftrightarrows\bb{H}(\cal{M}\msf{fld}_n,\cal{C}\msf{h}_\bb{Q}):\mrm{ev}_{\bb{R}^n}\]
\end{thm}

\begin{rmk}
One can replace $\cal{C}\msf{h}_\bb{Q}$ 
with a general symmetric monoidal $\infty$-category $\cal{V}$ 
as long as $\cal{V}$ is ``$\otimes$-presentable." 
For details, see \cite{AFtop}.
\end{rmk}

\subsection{Examples}

We compute factorization homology $\int_MA$ 
for simple choices of $M$ and $A$.

\begin{ex}
Take $M=\bb{R}^n$. 
Then $\cal{D}\msf{isk}_{n/\bb{R}^n}$ has a final object 
given by the identity map $\bb{R}^n=\bb{R}^n$. 
Thus the colimit is given by evaluation on $\bb{R}^n$,

\[\int_{\bb{R}^n}A=\msf{colim}\left(\cal{D}\msf{isk}_{n/\bb{R}^n}\rta\cal{D}\msf{isk}_{n}\xrta{A}\cal{V}\right)=A(\bb{R}^n)\]
\end{ex}

\begin{ex}
Take $M=\coprod\limits_I\bb{R}^n$. 
Then $\cal{D}\msf{isk}_{n/M}$ again has a final object 
and as above we obtain
\[\int_{\coprod\limits_I\bb{R}^n}A\simeq A(\coprod_I\bb{R}^n)\cong A(\bb{R}^n)^{\otimes I}\]
Here we are seeing the fact that 
\[\int_{(-)}A:\cal{M}\msf{fld}_n\rta\cal{C}\msf{h}_\bb{Q}\]
is a symmetric monoidal functor.
\end{ex}

\begin{ex}
Take $M=S^1$, as a framed manifold. 
Note that an $\cal{E}_1$-algebra $A$ is the same as an associative algebra 
$\bar{A}:=A(\bb{R}^1)$. 
We will use excision to compute $\int_{S^1}A$. 
Express $S^1$ as a collar-gluing 

\[S^1\cong\bb{R}\cup_{S^0\times\bb{R}}\bb{R}\]

By $\otimes$-excision, we have

\[\int_{S^1}A\simeq \left(\int_{\bb{R}}A\right)\bigotimes_{\left(\int_{S^0\times\bb{R}}A\right)}\left(\int_{\bb{R}}A\right)\simeq \bar{A}\bigotimes_{\bar{A}\otimes\bar{A}^\mrm{op}}\bar{A} \]
where we obtained 
$\bar{A}\otimes\bar{A}^\mrm{op}$ 
because the two copies of $\bb{R}^1$ in $S^0\times\bb{R}^1\subset S^1$ 
are oriented differently.
This is the \emph{Hochschild homology} of $A$. 

This appears as \cite[Thm. 3.19]{AFtop}.
\end{ex}

We have a functor 
$\msf{Alg}_n(\cal{V})\rta\cal{V}$ 
given by evaluating on $\bb{R}^n$. 
This is the forgetful functor.

\begin{defn}
The left adjoint to the forgetful functor is the free functor 
\[\bb{F}_n\colon\cal{V}\rta\msf{Alg}_n(\cal{V})\]
\end{defn}

\begin{ex}
Consider the free $n$-disk algebra on $V\in\cal{C}\msf{h}_\bb{Q}$. 
This sends a disjoint union $\coprod\limits_{k}\bb{R}^n$ to 
\[\bigoplus_{0\leq i}C_*(\msf{Emb}(\coprod_i\bb{R}^n,\coprod_k\bb{R}^n))\bigotimes_{\Sigma_i}V^{\otimes i}\]
\end{ex}

\begin{ex}
We can similarly define a free framed $n$-disk algebra. 
In the framed case,
\[\bb{F}_nV(\bb{R}^n)=\bigoplus_{i\geq 0}C_*\msf{Emb}^\mrm{fr}(\coprod_i\bb{R}^n,\bb{R}^n)\bigotimes_{\Sigma_i}V^{\otimes i}\]

Since $\msf{Emb}^\mrm{fr}(\coprod_i\bb{R}^n,\bb{R}^n)\simeq \msf{Conf}_i(\bb{R}^n)$, 
this agrees with the free $\cal{E}_n$-algebra on $V$.
\end{ex}

\begin{prop}
For $M$ a framed manifold, 
and $V\in\cal{C}\msf{h}_\bb{Q}$, 
we have
\[\int_M\bb{F}_nV\simeq\bigoplus_{0\leq i}C_*(\msf{Conf}_iM)\bigotimes_{\Sigma_i}V^{\otimes i}\]
\end{prop}

\noindent A similar statement is true in the non-framed case;
see \cite[Prop. 5.5]{AFtop}.

For $U\cong\coprod_I\bb{R}^n$ we have
\[\bb{F}_nV(U)=\left(\bigoplus_{i\geq 0}C_*(\msf{Conf}_i\bb{R}^n)\bigotimes_{\Sigma_i}V^{\otimes i}\right)^{\otimes I}\cong\bigoplus_{i\geq 0}C_*(\msf{Conf}_iU)\bigotimes_{\Sigma_i}V^{\otimes i}\]
Thus

\begin{align*}
\int_M\bb{F}_nV&=\underset{U\in\cal{D}\msf{isk}_{n/M}}{\msf{colimit}}\bigoplus_{i\geq 0}\left(C_*(\msf{Conf}_iU)\bigotimes_{\Sigma_i}V^{\otimes i}\right)\\
&=\bigoplus_{i\geq 0}\underset{U\in\cal{D}\msf{isk}_{n/M}}{\msf{colimit}}\left(C_*(\msf{Conf}_iU)\bigotimes_{\Sigma_i}V^{\otimes i}\right)
\end{align*}

Let $\msf{Disk}^\mrm{fr}_{n/M}$ 
denote the \emph{ordinary} category of framed $n$-disks in $M$. 
We will show things for the ordinary category $\msf{Disk}^\mrm{fr}_{n/M}$, 
instead of for the $\infty$-category $\cal{D}\msf{isk}^\mrm{fr}_{n/M}$. 
It turns out that this is sufficient
by \cite[Prop. 2.19]{AFtop}:

\begin{thm}
	The functor $\msf{Disk}^\mrm{fr}_{n/M}\rta \cal{D}\msf{isk}^\mrm{fr}_{n/M}$ is a localization. Hence factorization homology can be computed as a colimit over $\msf{Disk}^\mrm{fr}_{n/M}$.
\end{thm}

\noindent For a more direct proof in the $\infty$-category case, 
see \cite[Prop. 5.5.2.13]{HA}.

To compute this colimit, 
we use a hypercover argument. 
This is theorem A.3.1 in \cite{HA};
see also \cite{DI} or \cite[\S 4.2]{KnudsenConf}.

\begin{thm}[Seifert-van Kampen Theorem]
Let $X$ be a topological space. 
Let $\msf{Opens}(X)$ denote the poset of open subsets of $X$. 
Let $\cal{C}$ be a small category and let 
$F:\cal{C}\rta\msf{Opens}(X)$ be a functor. 
For every $x\in X$, let $\cal{C}_x$ 
denote the full subcategory of $\cal{C}$ 
spanned by those objects $C\in\cal{C}$ 
such that $x\in F(C)$. 
If for every $x\in X$, the simplicial set $N(\cal{C}_x)$ 
is weakly contractible, 
then the canonical map 
\[\underset{C\in\cal{C}}{\msf{colim}}\msf{Sing}(F(C))\rta\msf{Sing}(X)\]
exhibits the simplicial set 
$\msf{Sing}(X)$ as a homotopy colimit of the diagram 
$\{\msf{Sing}(F(C))\}_{C\in\cal{C}}$.
\end{thm} 

To use the Seifert-van Kampen theorem, consider the following commutative diagram,
\[\xymatrix{
\msf{Disk}^\mrm{fr}_{n/M}\arw[r]\arw[dr]\arw@/_2pc/[ddr]_{\msf{Conf}_i(-)} & \msf{Disk}^\mrm{fr}_{n}\arw[rr]^{\msf{Conf}_i(-)} & & \cal{C}\msf{h}_\bb{Q}\\
& \msf{Opens}(M)\arw[urr]_{\msf{Conf}_i(-)} & & \\
& \msf{Opens}(\msf{Conf}_iM)\arw@/_2pc/[uurr] & & 
}\]

Let $\bar{x}=(x_1,\dots,x_i)\in\msf{Conf}_iM$. 
The category $(\msf{Disk}^\mrm{fr}_{n/M})_{\bar{x}}$ 
contains embedded disks $U\hookrightarrow M$ 
so that $\{x_1,\dots,x_i\}$ is in $U$. 
By the Seifert-van Kampen theorem, if 
$B(\msf{Disk}^\mrm{fr}_{n/M})_{\bar{x}}\simeq *$, 
then 
\[\underset{U\in\msf{Disk}_{n/M}}{\msf{colimit}}\msf{Conf}_iU\simeq \msf{Conf}_iM\]

To show that this category is contractible, we will show it is cofiltered. 

\begin{defn}
A nonempty \emph{ordinary} category $\cal{C}$ is cofiltered if

1) for every pair $U,V\in\cal{C}$ there exists $W\in\cal{C}$ 
and maps $W\rta U$ and $W\rta V$, and

2) given two maps $u,v:X\rta Y$ in $\cal{C}$, 
there exists $Z\in \cal{C}$ and a map $w:Z\rta X$ so that $uw=vw$.
\end{defn}

\begin{proof}[Proof of Computation in the free case] 
Let $U,V\in(\msf{Disk}^\mrm{fr}_{n/M})_{\bar{x}}$. 
We need to find a finite disjoint union of euclidean spaces 
$W\rta M$ containing $(x_1,\dots,x_i)$ 
and maps $W\rta U$ and $W\rta V$. 
Note that $U\cap V$ contains $\bar{x}$, 
but may not be a disjoint union of euclidean spaces. 
However, we can find a small disk around 
each $x_i$ and still in $U\cap V$. 
The second condition is satisfied since 
$\msf{Disk}^\mrm{fr}_{n/M}$ is a poset.

Thus $(\msf{Disk}^\mrm{fr}_{n/M})_{\bar{x}}$ cofiltered, 
and hence contractible. 
Applying the Seifert-van Kampen theorem, 
(and adding in a few details about $V$) we get
\[\int_M\bb{F}_nV\simeq\bigoplus_{i\geq 0}C_*(\msf{Conf}_iM)\otimes_{\Sigma_i}V^{\otimes i}.\]
\end{proof}

For more details, see \cite[\S 5.2]{AFtop}.

\section{Exercises}

\subsubsection{Warm Up}

\begin{quest}\label{q-12.1}
What is the factorization homology 
\[\int_{[0,1]}A\]
of an associative algebra $A$?
\end{quest}

\begin{quest}\label{q-12.2}
What can you say about $\int_{S^2}A$ for $A$ a 2-disk algebra?
\end{quest}

\begin{quest}\label{q-3excision}
Make sense of the tensor product in excision.
How are things modules, algebras, and such?
\end{quest}

\subsubsection{Poincar\'e Duality for Factorization Homology}

Yesterday, you saw that $\Omega^nX$ was an $\bb{E}_n$-algebra.

\begin{quest}\label{q-12.4}
Show that $C_\bullet(\Omega^n X)$ is an $\bb{E}_n$-algebra.
Write $C_\bullet(\Omega^nX)$ and $\Omega^nX$ as $n$-disk algebras 
(valued in chain complexes and spaces, respectively).
Write them as factorization algebras as well.
\end{quest}

Note that 
\[\Omega^n X=\mathsf{Map}_c(\bb{R}^n, X).\]

\begin{quest}\label{q-12.5}
Let $X$ be an $n$-connective space.
Convince yourself that compactly supported maps 
\[\mathsf{Map}_c(-,X)\]
satisfies $\otimes$-excision. 
That is,
given a collar glueing 
\[M= U\bigcup_{V\times\bb{R}} U',\]
 there is an equivalence
\[\mathsf{Map}_c(U,X)\times_{\mathsf{Map}_c(V\times\bb{R},X)}\mathsf{Map}_c(U',X)\simeq\mathsf{Map}_c(M,X).\]
You can do this by just skimming the proof given in \cite[Lem. 4.5]{AFtop},
if you want; 
or by trying it out in a few easier cases.
\end{quest}

The following is a theorem of Salvatore \cite{Salvatore}, Segal \cite{SegalQFT}, and Lurie \cite[Thm. 5.5.6.6]{HA},
in various contexts. 
We are following the proof of Ayala-Francis \cite[Cor. 4.6]{AFtop}.

\begin{quest}[Nonabelian Poincar\'e Duality]\label{q-12.6}
Let $X$ be an $n$-connective space. 
Show that 
\[\int_M\Omega^nX\simeq\mathsf{Map}_c(M,X).\]
\end{quest}

Say $M$ is compact.
If $X$ an Eilenberg-MacLane space,
the right-hand side looks like cohomology. 
This motivates the relationship between nonabelian Poincar\'e duality
and usual Poincar\'e duality.

\subsubsection{Enveloping Algebras}

Let $\mathfrak{g}$ be a Lie algebra. 
Recall the Chevalley-Eilenberg complex 
$C_\bullet^\mathrm{Lie}(\mathfrak{g})$
from yesterday.

\begin{quest}\label{q-12.7}
Define a Lie algebra structure on 
\[\mathsf{Map}_c(\bb{R}^n,\mathfrak{g}).\]
Show that 
\[C_\bullet^\mathrm{Lie}(\mathsf{Map}_c(\bb{R}^n,\mathfrak{g}))\]
forms an $n$-disk algebra.
Call it $U_n\mathfrak{g}$.
\end{quest}

\begin{quest}\label{q-12.8}
Show that $U_1\mathfrak{g}$ is the enveloping algebra $U\mathfrak{g}$.
Check this with your understanding of $U\mathfrak{g}$
as an $\bb{E}_1$-algebra from yesterday.
\end{quest}

Thus $U_n\mathfrak{g}$ gives us a version of the enveloping algebra
in higher dimensions -a ``higher enveloping algebra"
This is a key example in field theory.
Many field theories have observables that look similarl to 
a higher enveloping algebra construction.
For example, any free theory has this property. 

\begin{quest}\label{q-3Ug}
Compute 
\[\int_M U_n\mathfrak{g}.\]
\end{quest}

\subsubsection{Spare Questions}

\begin{quest}
Show that $\int_MA$ has a canonical action of $\mathsf{Diff}(M)$.
\end{quest}

\begin{quest}
Let $H\subset \mathrm{GL}(n)$ be a sub-Lie-group.
You can think $SO(n)$ if you want.
Define a notion of an $H$-oriented TFT in the functorial setting.
\begin{enumerate}
\item Can you define a notion of an $H$-oriented $\bb{E}_n$-algebra?
\item How about an $H$-oriented factorization algebra?
\end{enumerate}
\end{quest}

\newpage
\part{Day 4}

\section{Functorial Things}

As we have seen in the past two days,
factorization algebras, $\bb{E}_n$-algebras, and $n$-disk algebras
are all related by the fundamental transformation 
of thinking of embedded disks, or points at their centers.
By similar reasoning, observables are sometimes called point observables.

To record what we learned,
we have a Costello-Gwilliam approach
and a special case when the field theory is topological.
The observables are a factorization algebra in general,
and in the topological case we 
get a $\bb{E}_n$-algebra in $\mathsf{Ch}$.

\begin{center}
\begin{tabularx}{0.8\textwidth} { 
  | >{\raggedright\arraybackslash}X 
  | >{\centering\arraybackslash}X 
  | >{\raggedleft\arraybackslash}X | }
 \hline
 & \begin{center}Costello-Gwilliam\end{center} & \begin{center}Topological Case\end{center} \\
 \hline
\begin{center}point observables\end{center}  & \begin{center}factorization algebra on spacetime\end{center}  & \begin{center}$\bb{E}_n$-algebra in $\mathsf{Ch}$ \end{center} \\
\hline
\end{tabularx}
\end{center}

Indeed, 
say our field theory is topological
and lives on $\bb{R}^n$.
The factorization algebra of observables only depends on the data of 
\[\mathcal{F}(\bb{D}^n)\in\mathsf{Ch}\]
and the multiplications coming from inclusions of disks into bigger disks.
If we draw this,
we see that $\mathcal{F}(\bb{D}^n)$ is the value at a disk centered at some point
and then inclusion of disks can be pulled out to be a bordism between
spheres.

\begin{center}
\includegraphics[scale=.4]{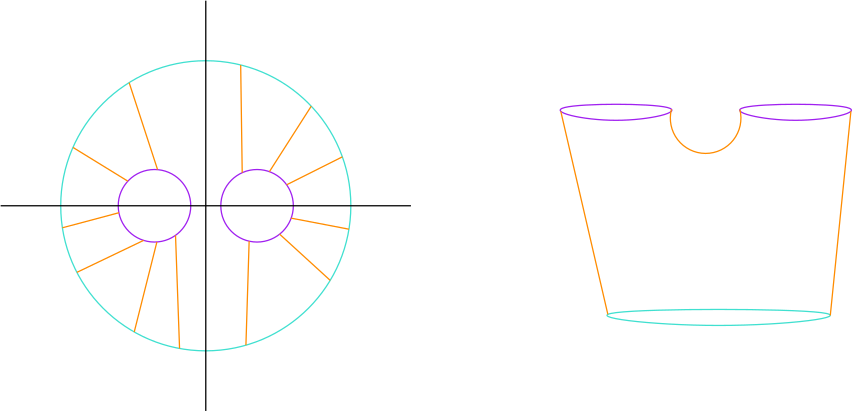}
\end{center}

This sphere is the linking sphere of the disk. 
The result, which I will denote
\[\mathcal{F}(\bb{D}^n)=\mathcal{Z}(S^{n-1})\]
for its dependence on the linking sphere
is the $\bb{E}_n$-algebra in $\mathsf{Ch}$.

\subsection{Definition of Line Operators}

Observables gave us some cool algebras to think about,
but there is some structure of the field theory that it misses.

\begin{ex}
Consider $G$-gauge theory on $M$. 
Given a loop $C$ in $M$,
we can define a map
\[\mathsf{Bun}^\nabla_GM\rta\bb{R}\]
sending a principal $G$-bundle $P\rta M$ 
to the trace of the holonomy map 
\[\mathrm{Hol}_C(P)\colon \mathfrak{g}\rta\mathfrak{g}.\]
This is called a Wilson loop operator.
\end{ex} 

So given a loop in spacetime, 
we get a function on fields,
which is something observables could know about,
but observables does not see any of the dependence of loops in $M$.

In general, 
this type of construction giving observables 
depending on loops or lines in spacetime
are called \emph{line operators}. 

\begin{rmk}
I know this is super vague.
It is my understanding that line operators
are still partially in the physics art stage rather than 
being fully mathematically understood.
\end{rmk}

Let's think about what type of structure
the set of line operators would have.

One piece of structure comes from stacking lines.
This gives us a way of composing line operators.
We should get a \emph{category} of line operators with

\begin{itemize}
\item objects: line operators (a pair of a line in spacetime and the operator), and
\item morphisms: $\mathsf{Hom}(L,L')$ is the set of point observables 
that can be inserted between $L$ and $L'$ to form a new line operator.
\end{itemize}

\begin{quest}
Is there any algebraic structure on the category of line operators?
\end{quest}

For observables,
the $\bb{E}_n$-algebra structure came from looking at the linking sphere
of the point we were at.
Let's replicate that with lines.
In $\bb{R}^n$, the linking sphere of a line is $S^{n-2}$. 
Now instead of including disks into bigger disks,
we have lines colliding 
(analogous to points colliding).
The result is a bordism between copies of $S^{n-2}$.
This is a $\bb{E}_{n-1}$-algebra structure.

\begin{claim}
For a topological theory on $\bb{R}^ n$,
line operators form an $\bb{E}_{n-1}$-monoidal category
$\mathcal{Z}(S^{n-2})$.
\end{claim}

\begin{center}
\begin{tabularx}{0.8\textwidth} { 
  | >{\raggedright\arraybackslash}X 
  | >{\centering\arraybackslash}X 
  | >{\raggedleft\arraybackslash}X | }
 \hline
 & \begin{center}Costello-Gwilliam\end{center} & \begin{center}Topological Case \end{center}\\
 \hline
\begin{center}point observables\end{center}   & \begin{center} factorization algebra on spacetime\end{center}   & \begin{center} $\mathcal{Z}(S^{n-1})$ an $\bb{E}_n$-algebra in $\mathsf{Ch}$\end{center}   \\
 \hline
\begin{center}line operators\end{center}   & \begin{center}?\end{center} & \begin{center}$\mathcal{Z}(S^{n-2})$ an $\bb{E}_{n-1}$-algebra in $\mathsf{Cat}$\end{center}  \\
\hline

\end{tabularx}
\end{center}

\subsection{Line Operator Constructions}

How do we say anything about line operators on the non-topological side?

\begin{idea}
Let $F$ be a field theory on $M$.
There is an ansatz from physics that 
a line operator for $C\subset M$
is the same as the data of boundary theory for 
$F$ restricted to $M\setminus \mathsf{Norm}(C)$.
\end{idea}
That is,
the field theory works the same away from $C$ 
and has a ``defect" at $C$;
see \cite{CEG} for a discussion in this vein. 

Thinking in terms of observables, 
in the topological $\bb{R}^n$ case,
we would like an $\bb{E}_n$-algebra 
away from $C$. 

\begin{quest}
What is the data of observables of a boundary theory for 
$F$ on $M\setminus \mathsf{Norm}(C)$?
\end{quest}

The following is part of \cite[Prop. 4.8]{AFT2}.

\begin{thm}[Ayala-Francis-Tanaka]
Let $A$ be an $\bb{E}_n$-algebra.
The data needed to produce a new $\bb{E}_n$-algebra 
that agrees with $A$ on $\bb{R}^n\setminus \bb{R}^k$ is 
an object in
\[\mathsf{LMod}\left(\int_{S^{n-k-1}\times\bb{R}^{k+1}}A\right).\]
\end{thm}

\begin{rmk}
To state this theorem,
we are using the $\bb{E}_1$-algebra structure on 
\[\int_{S^{n-k-1}\times{R}^{k+1}}A.\]
This comes from the $\bb{R}^{k+1}$ direction by stacking.
This is related to the exercise from yesterday 
on making sense of excision.
\end{rmk}

For line operators, 
we are looking at modules over 
\[\int_{S^{n-2}}A.\]
This $S^{n-2}$ is the linking sphere of the line $C$ that we encountered before.

Thus, we have an approximation to the category of line operators
\[\mathsf{LineOp}\approx \mathsf{Mod}\left(\int_{S^{n-2}}A\right).\]

\begin{rmk}
This is not a good approximation for non-perturbative field theories.
\end{rmk}

\section{Atiyah-Segal Approach to TFTs}
Returning to the topological side,
we start to see a pattern.
Continuing down,
we could ask for the date of 
an $\bb{E}_n$-algebra $\mathcal{Z}(S^{n-1})$,
an $\bb{E}_{n-1}$-monoidal category $\mathcal{Z}(S^{n-2})$,
and so on.

We will use a functorial TQFT perspective to 
package together into a single functor.

\begin{notation}
For $X$ an oriented manifold, let $\overline{X}$ denote $X$ with the opposite orientation. 
\end{notation}

The following is \cite[Def. 1.1.1]{LurieCobordism}.

\begin{defn}
Let $n$ be a positive integer. We define a category $\mathbf{Cob}(n)$ as follows:
\begin{itemize}
\item objects: $(n-1)$-dimensional oriented manifolds 
\item morphisms from $M$ to $N$ are given by equivalence classes of $n$-dimensional oriented manifolds with boundary $B$ together with an an orientation-preserving diffeomorphism 
\[\del B\simeq \overline{M}\sqcup N\] 
Two morphisms $B$ and $B'$ are equivalent if there is an orientation-preserving diffeomorphism $B\rta B'$ that restricts to the identity on the boundaries,
\[\xymatrix{
B\arw[r]\arw[d]^\simeq  & B'\arw[d]^\simeq \\
\overline{M}\sqcup N\ar@{=}[r] & \overline{M}\sqcup N
}\]
\end{itemize}

Composition is given by gluing cobordisms along their shared boundary.
View $\mathbf{Cob}(n)$ as a symmetric monoidal category under disjoint union.
\end{defn}

For $k$ a field, let $\mathbf{Vect}(k)$ denote the symmetric monoidal category of $k$-vector spaces with tensor product.

The following comes from \cite{Atiyah,Segal}.

\begin{defn}[Atiyah, Segal]
Let $k$ be a field. A \emph{topological field theory} (TFT) of dimension $n$ is a symmetric monoidal functor $Z\colon\mathbf{Cob}(n)\rta\mathbf{Vect}(k)$.
\end{defn}

In particular, $Z(\emptyset)=k$.

\begin{lem}
The value $Z(S^{n-1})$ is an $\bb{E}_n$-algebra.
\end{lem}

\begin{proof}
The pair of pants bordism gives the 
\[\bb{E}_n(2)\otimes Z(S^{n-1})^{\otimes 2}\rta Z(S^{n-1})\]
map.
More legs, means more points.
\end{proof}

\begin{rmk}Note that Atiyah's definition considers all possible spacetimes at once, instead of working on one specific $n$-manifold at a time. We could instead consider the category of bordism submanifolds within a fixed manifold.
 \end{rmk}

\begin{notation}
For $V$ a $k$-vector space, let $V^\vee$ denote the linear dual, $V^\vee\colon=\mrm{Hom}(V,k)$.
\end{notation}

Let $Z$ be an $n$-dimensional TFT. Given an oriented $(n-1)$-manifold $M$, the product manifold $M\times[0,1]$ can be viewed as a morphism in $\mathbf{Cob}(n)$ in multiple ways.

\begin{itemize}
\item As a morphism from $M\rta M$, the product $M\times[0,1]$ maps to the identity map \[\mrm{id}\colon Z(M)\rta Z(M).\]

\item As a morphism $M\sqcup\overline{M}\rta\emptyset$, the product $M\times[0,1]$ determines an evaluation map 
\[\mrm{ev}\colon Z(M)\otimes Z(\overline{M})\rta k.\]

\item As a morphism $\emptyset\rta \overline{M}\sqcup M$, the product $M\times[0,1]$ determines a coevaluation map
\[\mrm{coev}\colon k\rta Z(\overline{M})\otimes Z(M).\]
\end{itemize}

Recall that a pairing $V\otimes W\rta k$ is perfect if it induces an isomorphism $V\rta W^\vee$.
The following is \cite[Prop. 1.1.8]{LurieCobordism}.

\begin{prop}\label{dual}
Let $Z$ be a topological field theory of dimension $n$. For every $(n-1)$-manifold $M$, the vector space $Z(M)$ is finite dimensional. The evaluation map $Z(M)\otimes Z(\overline{M})\rta k$, induced from the cobordism $M\times[0,1]$, is a perfect pairing.
\end{prop}

\section{Classifying Topological Field Theories}

\subsection{Low Dimensions}

\begin{ex}[Dimension 1]
Let $Z\colon\mathbf{Cob}(1)\rta\mathbf{Vect}(k)$ be a 1-dimensional TFT. Let $P$ denote a single point with positive orientation and $Q=\overline{P}$. Let $Z(P)=V$. This finite-dimensional vector space, determines $Z$ on objects. By Proposition \ref{dual}, $Z(Q)=Z(\overline{P})=V^\vee$. A general object of $\mathbf{Cob}(1)$ looks like
\[M=\coprod_{S_+}P\sqcup \coprod_{S_-}Q\]
for $S_+,S_-$ sets. Since $Z$ is symmetric monoidal, we have
\[Z(M)=\coprod_{S_+}V\sqcup \coprod_{S_-}V^\vee\]
What about morphisms? A morphism in $\mathbf{Cob}(1)$ is a 1-dimensional manifold with boundary $B$. Using the monoidal structure, it suffices to describe $Z(B)$ where $B$ is connected. There are five possibilities
\begin{itemize}
\item $B$ is an interval viewed as a morphism $P\rta P$. Then $Z(B)=\mrm{Id}_V$.
\item $B$ is an interval viewed as a morphism $Q\rta Q$. Then $Z(B)=\mrm{Id}_{V^\vee}$.
\item $B$ is an interval viewed as a morphism $P\sqcup Q\rta \emptyset$. Then $Z(B)\colon V\otimes V^\vee\rta k$. By Proposition \ref{dual}, this is the canonical pairing of $V$ and $V^\vee$,
\[Z(B)(x\otimes f)=f(x)\]
Under the isomorphism $V\otimes V^\vee\cong \mrm{End}(V)$, the morphism $Z(B)$ corresponds to taking the trace.
\item $B$ is an interval viewed as a morphism $\emptyset\rta P\sqcup Q$. Then $Z(B)$ is the map $$k\rta V\otimes V^\vee\cong\mrm{End}(V)$$ 
sending $\lambda\in k$ to $\lambda\mrm{Id}_V$.
\item $B=S^1$ is a circle viewed as a morphism $\emptyset\rta\emptyset$. Then $Z(S^1)$ is a linear map $k\rta k$; i.e., multiplication by some $\gamma\in k$. To determine $\gamma$, view $S^1$ as the union of two semi-circles along $P\sqcup Q$. This determines a decomposition of $S^1$ into the composite of two cobordisms, 
\[\emptyset\rta P\sqcup Q\rta\emptyset\]
By the above cases, $Z$ maps this to the composite
\[k\rta \mrm{End}(V)\rta k\]
of the map $\lambda\mapsto \lambda\mrm{Id}_V$ and the trace map. Thus $Z(S^1)$ is the scaling by $\mrm{Tr}(\mrm{Id}_V)=\dim(V)$ map.
\end{itemize}
\end{ex}
\begin{rmk}
In 1-dimension, the observables are given by $Z(S^0)=\mrm{End}(V)$. Notice that this has the structure of an associative algebra.
\end{rmk}
Thus in dimension 1, we see that the vector space $Z(P)$ determines the TFT $Z$. Does every $V\in\mathbf{Vect}(k)$ appear as $Z(P)$ for some 1-dimensional TFT? Nope, only the finite-dimensional ones. We get an equivalence of categories
\[\mrm{Fun}^\otimes(\mathbf{Cob}(1),\mathbf{Vect}(k))\rta\mathbf{Vect}^\mrm{fin}(k)\]
by evaluating on the point.

Let's try to do something similar in dimension 2.

\begin{ex}[Dimension 2]
Following \cite[Ex. 1.1.11]{LurieCobordism} we have the following.
Let $Z$ be a 2-dimensional TFT. The only objects in $\mathbf{Cob}(2)$ are the empty set and disjoint unions of copies of $S^1$. We do not get a new object $\overline{S^1}$ since the circle has an orientation-reversing diffeomorphism. The observables, $A=Z(S^1)$ determines $Z$ on objects.  

What about morphisms? A morphism in $\mathbf{Cob}(2)$ is 2-dimensional oriented manifold with boundary. 
\begin{itemize}
\item The pair of pants cobordism determines a map $m\colon A\otimes A\rta k$. One can check that $m$ defines a commutative, associative multiplication on $A$. 
\item The disk $\bb{D}^2$ viewed as a cobordism $S^1\rta \emptyset$ determines a linear map $\mrm{Tr}\colon A\rta k$.
\item The disk $\bb{D}^2$ viewed as a cobordism $\emptyset\rta S^1$ determines a linear map $k\rta A$. The image of $1\in k$ under this map acts as a unit for the multiplication. Indeed, we can glue $\bb{D}^2$ to one of the legs of the pants. The resulting manifold is diffeomorphic to $S^1\times[0,1]$. But $S^1\times[0,1]$ maps to $\mrm{Id}_A$ under $Z$.
\end{itemize}
Note that the composite of 
\[A\otimes A\xrta{m} A\xrta{\mrm{Tr}}k\]
comes from the cobordism $S^1\times[0,1]$ viewed as a map $S^1\sqcup S^1\rta \emptyset$. By Proposition \ref{dual}, the map $\mrm{Tr}\circ m$ is a nondegenerate pairing.
\end{ex}
\begin{defn}
A \emph{commutative Frobenius algebra} over $k$ is a finite-dimensional commutative $k$-algebra $A$, together with a linear map $\mrm{Tr}\colon A\rta k$ such that the bilinear form $(a,b)\mapsto\mrm{Tr}(a,b)$ is nondegenerate.
\end{defn}

\begin{thm}[Folklore]
The category of 2-dimensional oriented TFTs is equivalent to the category of Frobenius algebras.
\end{thm}

A detailed proof can be found in \cite{Abrams} and a good expository account can be found in \cite{Kock}.
See \cite{SP} for details in the fully extended case.

\begin{rmk}
In 2-dimensions, the observables $Z(S^1)$ of a TFT $Z$ has the structure of a Frobenius algebra.
\end{rmk}

\section{Exercises}

\begin{quest}\label{q-16.1}
Let $V$ be a finite dimensional vector space.
Let $Z_V$ be the 1-dimensional TQFT determined by $V$.
What is the corresponding associative algebra of observables,
in terms of $V$?
\end{quest}

\begin{quest}\label{q-16.2}
Let $Z$ be an $n$-dimensional field theory.
Show that $Z(S^{n-1})$ acts on $Z(N)$ for every $n$-manifold $N$.
\end{quest}

\begin{quest}\label{q-16.3}
For $G$-gauge theory on $M$,
we defined a Wilson loop operator.
Given a representation $\rho$ of $G$,
define a version of a Wilson loop operator on $\mathsf{Gauge}_GM$.
This should be an assignment of a real number to a principal $G$-bundle
$P\rta M$ with connection.
As a first step, you should build the associated vector bundle to $P$ using $\rho$.
\end{quest}

\begin{quest}\label{q-16.5}
Let $A$ be an $\bb{E}_n$-algebra.
For example, $A=\mathsf{Obs}^q$ the observables of a TQFT on $\bb{R}^n$.
The line operators of this field theory should 
be a $\bb{E}_{n-1}$-monoidal category. 
As guesses, 
build two different categories
out of $A$.
Do these categories have any monoidal structure?
\end{quest}

\subsubsection{Classification of 3d TQFTs}

We saw that, 
up to reversing orientation and taking disjoint unions,
\begin{itemize}
\item[] $\mathsf{Cob}(1)$ had a unique object $P$, the point, and
\item[] $\mathsf{Cob}(1)$ had a unique object $S^1$.
\end{itemize}

Examples of objects in $\mathsf{Cob}(3)$ include the torus $T^2$ and the sphere $S^2$.

\begin{quest}\label{q-16.6}
Using the classification of surfaces, 
describe the set of objects of $\mathsf{Cob}(3)$.
\end{quest}

\begin{quest}\label{q-16.7}
Recall the operation 
of connected sum for oriented closed 2-manifolds.
Show that under connect sum, one can construct every object of $\mathsf{Cob}(3)$
from $T^2$ and $S^2$.
\end{quest}

Note that objects of $\mathsf{Cob}(3)$ 
determine morphisms in $\mathsf{Cob}(2)$;
we can view a closed $2$-manifold as a
cobordism from the emptyset to itself.

\begin{quest}\label{q-16.8}
As morphisms in $\mathsf{Cob}(2)$,
what does the operation you thought of in
Question \ref{q-16.7} correspond to categorically?
\end{quest}

This breaking complicated 2-dimensional manifolds
down into easy pieces (like $S^2$ and $T^2)$ 
is very helpful. 
If we want a classification of 3-dimensional TQFTs,
we would like to be able to do this in $\mathsf{Cob}(3)$.

\begin{quest}\label{q-16.9}
Can you think of a way to encode the operation
from Question \ref{q-16.7} into a categorical setting 
involving 3-dimensional bordisms?
(Any guesses or ideas are fine, this one is to just get you thinking!)
\end{quest}

\newpage
\part{Day 5}

Yesterday we talked about how for TQFTs
the factorization algebra of observables 
(which is $\bb{E}_n)$
can be situated in the Atiyah-Segal definition of a TQFT.
We ended by classifying TQFTs in dimensions 1 and 2.

\section{Cobordism Hypothesis - Higher Dimensions}
The problem when we try to classify
 TFTs in higher dimensions is that the 
 objects become too complicated. 
 Up to reversing orientation and taking disjoint unions, 
 the categories $\mathbf{Cob}(1)$ and $\mathbf{Cob}(2)$
  have a unique object, $P$ and $S^1$, respectively. 
  For $n=3$, there are infinitely many oriented 2-manifolds, 
  one for each genus $g$. 
  We do not think of genus $g$ surfaces as being that complicated. 
  In fact, we usually think of $\Sigma_g$, the genus $g$ surface, 
  as coming from $g$ connect sums of the torus. 
  A closely related way to say this,
   is that $\Sigma_g$ has a relatively easy handle-body decomposition. 
   But what happens when we view $\Sigma_g$ under its handle-body decomposition? 
   We are really viewing it as a composition of cobordisms; 
   i.e., as a morphisms in $\mathbf{Cob}(2)$. 
   Similarly, when we tried to understand the value of a 2-dimensional field theory on $S^1$, 
   we broke $S^1$ into the union of two semi-circles, 
   that is to say, into its handle-body decomposition.

If we want to be understand an 
$n$-dimensional field theory by breaking manifolds down, 
using their handle-body decompositions, 
into lower-dimensional manifolds, 
we need the TFT to know about manifolds of dimension $<n-1$. 
In particular, we would like some sort of data assigned to 
every $(n-2)$-dimensional manifold and 
we would like this data to have something to do with the values on 
$(n-1)$-manifolds. 
In particular, from our discussion yesterday,
we expect $S^{n-2}$ to be assigned a category. 

The way to encode all this data is the language of higher categories.

\begin{defn}
A \emph{strict $n$-category} is a category $\cal{C}$ enriched over $(n-1)$-categories.
\end{defn}

For $n=2$, this means that for objects $A,B\in\cal{C}$ the morphisms $\msf{Hom}_\cal{C}(A,B)$ is itself a category.

\begin{ex}
The strict $2$-category $\mathbf{Vect}_2(k)$ has objects cocomplete $k$-linear categories and morphisms 
\[\mrm{Hom}_{\mathbf{Vect}_2(k)}(C,D)=\mathbf{Fun}_k^{\mrm{cocon}}(C,D)\]
the functor category of cocontinuos, $k$-linear, functors.
\end{ex}

\begin{ex}
The strict 2-category $\mathbf{Cob}_2(n)$ has 
\begin{itemize}
\item objects: closed, oriented manifolds of dimension $n-2$.
\item morphisms: $\mrm{Hom}_{\mathbf{Cob}_2(n)}(X,Y)=:\cal{C}$ should be the category with
\begin{itemize}
\item objects: cobordisms $X\rta Y$
\item morphisms: $\mrm{Hom}_\cal{C}(B,B')$ is equivalence classes of bordisms $X$ from $B\rta B'$ 
\end{itemize}

\end{itemize}

The big problem here is making the composition law strictly associative. One would like to define composition by gluing bordisms, but get messed up in defining a smooth structure on the result, and things that used to be equalities are now just homeomorphisms. The solution will be to get rid of the ``strictness" and move to $(\infty,2)$-categories.
\end{ex}

Let $\mathsf{Cob}_n(n)$ denote an $(\infty,n)$-category version of the cobordism category, 
see Calaquee-Scheimbauer \cite{Calaque-Scheimbauer}.

\begin{defn}
Let $\cal{C}$ be a symmetric monoidal $(\infty,n)$-category. An \emph{extended} $\cal{C}$-valued $n$-dimensional TFT is a symmetric monoidal functor
\[Z\colon\msf{Cob}_n(n)\rta\cal{C}\]
\end{defn}

\begin{ex}
Take $\mathcal{C}$ so that 
$\Omega^n\mathcal{C}=\bb{C}$,
$\Omega^{n-1}\mathcal{C}=\mathsf{Ch}_\bb{C}$,
and
$\Omega^{n-2}\mathcal{C}=\mathsf{LinCat}_\bb{C}$.
Then $Z(S^{n-2})$ is an $(\infty,1)$-category.
It has an $\bb{E}_{n-1}$-monoidal structure from the pair of pants bordism.
This is \emph{line operators}.
Similarly, $Z(S^{n-k})$ is an $\bb{E}_{n-k-1}$-monoidal $(\infty,n-k-1)$-category.
It describes $(k+1)$-dimensional defects (i.e. operators).
\end{ex}

The purpose of this definition is to allow us to reduced $n$-dimensional TFTs down to information about 1-dimensional TFTs. As we saw before, a 1-dimensional TFT is determined by its value on a point. Thus we might make the following guess.

\begin{guess}
An extended field theory is determined by its value on a single point. Moreover, evaluation on a point determines an equivalence of categories between TFTs valued in $\cal{C}$ and $\cal{C}$.
\end{guess}

There are two problems with this guess.

\begin{enumerate}
\item Even in 1-dimension, not every vector space determined a TFT. 
We needed to restrict to finite-dimensional ones. 
The analogue in higher dimensions will be something called ``fully dualizable objects."
\item Orientation in dimension 1 is the same as a framing. 
This is not true in higher dimensions. 
We actually wanted a framing, 
not just an orientation so that we could say that locally $M^k$ was canonically diffeomorphic to $\bb{R}^k$ 
(via the exponential map). 
Thus we need a version of $\msf{Cob}_n(n)$ that works with framed manifolds instead of oriented ones.
\end{enumerate}

The following conjecture is due to Baez and Dolan \cite{Baez-Dolan}.

\begin{conj}[Cobordism Hypothesis: Framed Version]
Let $\cal{C}$ be a symmetric monoidal $(\infty,n)$-category with duals. Then the evaluation functor $Z\mapsto Z(*)$ induces an equivalence
\[\msf{Fun}^\otimes(\mathbf{Bord}_n^\mrm{fr},\cal{C})\rta\cal{C}^\sim\]
between framed extended $n$-dimensional TFTs valued in $\cal{C}$ and the fully dualizable subcategory of $\cal{C}$. 
\end{conj}

Partial and complete proofs are due to Hopkins-Lurie, Lurie \cite{LurieCobordism}, Grady-Pavlov \cite{Grady-Pavlov}.
Many others have worked on this as well.
See \cite{FreedCob} for an expository account of the cobordism hypothesis
and quantum field theory.

\begin{rmk}
This leads to the natural question,
given $Z(\mathrm{pt})$, 
how does one obtain $Z(N)$?
It is expected that there is a version of factorization homology 
for $(\infty,n)$-categories so that 
\[\int_N Z(\mathrm{pt})=Z(N)\]
for all manifolds $N$ of dimension $0,\dots,n$. 
This is discussed in work of Ayala-Francis;
see \cite{AFcob},
where it is shown how to prove the cobordism hypothesis assuming the 
existence of  an upgraded version of factorization homology
known as ``$\beta$-factorization homology."
\end{rmk}

Take $\mathcal{C}$ to be a suitable choice of an $(\infty,n)$-category
of algebras up to Morita equivalences
\[\mathcal{C}=\mathsf{Morita}_n.\]

See Scheimbaurer's thesis \cite{Sche} for a factorization homology
reconstruction of a fully extended TQFT with Morita target.

The fully dualizable objects of $\mathsf{Morita}_{n-1}$
are certain types of $\bb{E}_{n-1}$-algebras.
Thus, you can describe an $n$-dimensional TQFT
by just giving an $\bb{E}_{n-1}$-algebra (satisfying certain conditions)
that the field theory assigns to a point.

Usually, we are thinking of an $n$-dimensional TQFT
as corresponding to its $\bb{E}_n$-algebra of observables,
so what is up with the $\bb{E}_{n-1}$-algebra?
How do we get from $Z(\mathrm{pt})$ to $Z(S^{n-1})$?

\subsection{Drinfeld Centers}

To answer this question,
we are going to use a version of 
a notion Delaney talked about yesterday \cite{Delaney}.

You saw yesterday that Drinfeld centers created braided monoidal structures from just monoidal structures.
You can think of that as saying that the center
of an $\bb{E}_1$-category is $\bb{E}_2$.
More generally,
we have the higher version of the Deligne conjecture due to 
Kontsevich, \cite{Ko}. 

\begin{thm}[Deligne Conjecture; Lurie]
The $\bb{E}_n$-Drinfeld center of an $\bb{E}_n$-category is and $\bb{E}_{n+1}$-category.
\end{thm}
\noindent This was proven in full generality in \cite[\S 5.3]{HA}. 
Specifically, see \cite[Cor. 5.3.1.15]{HA}.

One can show that $\bb{E}_n$-Drinfeld center of 
$Z(\mathrm{pt})$ is $Z(S^{n-1})$,
answering our previous question.

As a good geometric example, we have the following theorem.

\begin{thm}[Ben-Zvi, Francis, Nadler]
Let $X$ be a perfect stack.
The center of quasi-coherent sheaves on $X$ is 
sheaves on the free loop space,
\[\mathsf{Cent}_{\bb{E}_n}(QC(X))\simeq QC(\mathcal{L}^nX).\]
\end{thm}

See \cite[Thm. 1.7 and Cor. 5.12]{BZFN}.

\section{Holomorphic Field Theories}

Now I want to switch gears and talk about 
non-topological field theories.
These will provide our first example of factorization algebras
that are not $\bb{E}_n$-algebras.

Recall that a function 
\[f\colon\bb{C}\rta\bb{C}\]
is \emph{holomorphic}
if it is complex differentiable at every point.
Equivalently,
if we write
\[f(x+iy)=u(x,y)+iv(x,y),\]
and assume $f$ is continuous,
then $f$ is holomorphic if and only if $f$ 
satisfies the Cauchy-Riemann equations
\[\frac{\del u}{\del x}=\frac{\del v}{\del y}\]
and 
\[\frac{\del u}{\del y}=-\frac{\del v}{\del x}.\]
This is the Looman-Menchoff theorem.

Basically these are incredibly nice complex functions.
In order to make sense of holomorphic conditions,
our holomorphic field theories will have spacetime $\bb{C}^n$.

Topological field theories were particularly nice
because of their invariance property on observables,
\[\mathsf{Obs}^q(D_1)\xrta{\sim}\mathsf{Obs}^q(D_2).\]
To be able to get somewhere with holomorphic theories,
we will additionally assume an invariance property:
translation invariance.

Let $V$ be a holomorphically translation-invariant vector bundle on $\bb{C}^n$.
This means that we are given a holomorphic isomorphism between $V$ 
and a trivial bundle.
Our space of fields will be 
\[\Omega^{0,*}(\bb{C}^n,V).\]
Note that by Dolbeault's theorem,
this complex has cohomology
\[H^{ 0,*}(\bb{C}^n,V)=H^*(\bb{C}^n,\Omega^0\otimes V)\]
where $\Omega^0$ is the sheaf of holomorphic functions on $\bb{C}^n$. 
Thus, our space of fields is a derived model for the mapping space
of holomorphic functions
\[\mathsf{Map}_\mathrm{hol}(\bb{C}^n,V).\]

Let 
\[\eta_i=\frac{\del}{\del\bar{z}_i}\vee (-)\colon\Omega^{0,k}(\bb{C}^n,V)\rta\Omega^{0,k-1}(\bb{C}^n,V)\]
be the contraction operator.

\begin{defn}
A field theory on $\Omega^{0,*}(\bb{C}^n,V)$ 
is \emph{holomorphically translation invariant} if 
the action functional 
\[S\colon \Omega^{0,*}(\bb{C}^n,V)\rta \bb{C}\]
is translation-invariant and satisfies 
\[\eta_i S=0\]
for all $i=1,\dots,n$.
\end{defn}

As with topological theories,
we are interested in how the 
holomorphically translation-invariant 
condition on a field theory 
impacts the factorization algebra of observables.

To make the following definition precise
requires more background and time than we have.
Here is the idea;
details can be found in \cite[Def. 8.1.1]{CG1} and \cite[Ch 5. Def. 1.1.1]{CG1}.

\begin{defn}
A factorization algebra $\mathcal{F}$ on $\bb{C}^n$
is \emph{holomorphically translation invariant}
if we have isomorphisms
\[\tau_x\colon\mathcal{F}(U)\simeq \mathcal{F}(\tau_x U)\]
for all $x\in\bb{R}^n$ and open $U\subset\bb{R}^n$.
These isomorphisms are required to vary
holomorphically in $x$ and satisfy 
\[\tau_x\circ \tau_y=\tau_{x+y}\]
and commute with the factorization algebra maps.
\end{defn}

Note that $\mathsf{Obs}^q$ will be a 
factorization algebra on $\bb{C}^n$.
The following is \cite[Prop. 9.1.1.2]{CG2}.

\begin{thm}[Costello-Gwilliam]
The observables $\mathsf{Obs}^q$ 
of a holomorphically translation-invariant 
field theory
is a holomorphically translation-invariant factorization algebra.
\end{thm}

Say we are working over $\bb{C}$.
Assume that our field theory additionally
is $S^1$-invariant.
That is,
that there is an $S^1$ action on $V$
which,
together with the $S^1$ action on $\Omega^{0,*}(\bb{C})$,
gives an action on the space of fields,
and all the structures of the field theory are invariant under this.

In this situation, we will get a nice algebraic 
description of the observables,
like we did for locally constant factorization algebras
as $\bb{E}_n$-algebras. 

The following is \cite[Ch. 5, Thm. 0.1.3]{CG1}.

\begin{thm}[Costello-Gwilliam]\label{thm-CG6}
A holomorphically translation-invariant and $S^1$-invariant 
factorization algebra on $\bb{C}$ 
determines what is called a vertex algebra.
\end{thm}

\begin{defn}[Borcherds]
A \emph{vertex algebra} is the following data:
\begin{itemize}
\item a vector space $V$ over $\bb{C}$ (the state space);
\item a nonzero vector $|\Omega\rangle \in V$ (the vacuum vector);
\item a linear map $T\colon V\rta V$ (the shift operator);
\item a linear map 
\[Y(-,z)\colon V\rta\mathsf{End}V[[z,z^{-1}]]\]
\end{itemize}
such that 
\begin{itemize}
\item (vacuum axiom) $Y(|\Omega\rangle ,z)=\mathrm{id}_V$ and 
\[Y(v,z)|\Omega\rangle \in v+zV[[z]]\]
for all $v\in V$;
\item (translation axiom) $[T,Y(v,z)]=\del_zY(v,z)$ for every $v\in V$ 
and $T|0>=v$;
\item (locality axiom) for any pair of vectors $v,v'\in V$,
there exists a nonnegative integer $N$ such that 
\[(z-w)^N[Y(v,z),Y(v',w)]=0\]
as an element of $\mathsf{End}V[[z^{\pm 1},w^{\pm 1}]]$.
\end{itemize}
\end{defn}

See \cite{BenZviFrenkel} for a good reference on vertex algebras.

\begin{rmk}
Vertex algebras were around before factorization algebras.
They have the benefit of being super computational.
If you are really good at power series manipulations,
or more familiar with representation theory methods,
vertex algebras might be better suited for you.
Factorization algebras are more geometric
and closer to the topologists $\bb{E}_n$-algebras.
\end{rmk}

As motivation for the proof, 
recall how we got an $\bb{E}_n$-algebra $A$
from a locally constant factorization algebra $\mathcal{F}_A$ on $\bb{R}^n$. 
The underlying space of $A$ is given by 
\[\mathcal{F}_A(U)\]
where $U$ is any disk in $\bb{R}^n$. 
For convenience, 
we can take $U=B_1(0)$,
the unit ball around the origin.
The structure maps are from the inclusions 
of disjoint disks in to $B_1(0)$,
\[\mathcal{F}_A(B_r(0))\otimes \mathcal{F}_A(B_{r'}(0))\rta\mathcal{F}_A(B_1(0)).\]

We will get a vertex algebra from
a holomorphically translation-invariant 
factorization algebra by a similar process,
with one important difference.

\begin{center}
Vertex algebras act more like Lie algebras.
\end{center}

They have structures like Lie brackets,
rather than multiplications like groups do.
If evaluating on $B_1(0)$
gave us a group like structure,
than to get a Lie algebra structure,
we should somehow take the tangent space at 0.
This is in analogy with how given
a Lie group $G$,
one obtains the Lie algebra
as $T_eG$.

\begin{proof}[Proof Idea]
Let $\mathcal{F}$ be the 
holomorphically translation-invariant and $S^1$-invariant
factorization algebra on $\bb{C}$.
Let
\[\mathcal{F}_k(B_r(0))\]
be the weight $k$ eigenspace of the $S^1$ action.
To zoom in on $0$,
we take the limit
\[V_k=\mathrm{lim}_{r\rta 0} H^\bullet(\mathcal{F}_k(B_r(0)).\]
(CG assume the maps in this limit are quasi-isomorphisms.)
(They also assume that $V_k=0$ for $k>>$.)

The underlying vector space of the
vertex algebra will be 
\[V=\bigoplus_{k\in\bb{Z}}V_k.\]

The vacuum element is given by the image of 
the unit in $\mathcal{F}(\emptyset)$.

The translation map is given by 
the derivation $\frac{\del}{\del z}$ 
from infinitesimal translation in the $z$ direction. 
The state-field map 
\[Y\colon V\rta\mathsf{End}(V)[[z,z^{-1}]]\]
comes from the inclusion of two disjoint disks into a bigger disk,
and expanding the holomorphic map into a Laurent series.
\end{proof}

\begin{rmk}
One can attach a factorization algebra to a certain type of 
vertex algebra;
see \cite{Brue}.
\end{rmk}

Factorization homology in the language of vertex algebras
is related to ``conformal blocks."

\subsection{Examples}

\begin{ex}
A commutative ring $V$ with derivation $T$ 
determines a vertex algebra 
with state space $V$,
translation operator $T$,
and state-field correspondence 
\[Y(u,z)v=uv.\]
\end{ex}

\begin{ex}
An important example of a vertex algebra 
is the Virasoro vertex algebra.
For a detailed description of the associated factorization algebra
and a great example of Theorem \ref{thm-CG6},
see \cite{BrianVir}.
\end{ex}

\begin{ex}
Given a manifold $X$,
differential operators $\mathsf{Diff}_X$ is an associative algebra,
so a type of factorization algebra on $\bb{R}$.
A 2-dimensional analogue would 
be something like a factorization algebra on $\bb{R}^2)$;
for example, a vertex algebra.
Such a vertex algebra was constructed by 
Malikov-Schechtman-Vaintrob \cite{chiral}
and is called \emph{chiral differential operators}.
The factorization algebra analogue was build 
by Gorbounov-Gwilliam-Williams \cite{GGW}. 
This is an example of a factorization algebra
built locally on the target from a Lie algebra in an enveloping algebra-esque construction.
\end{ex}

\subsection{Stolz-Teichner Program}

The field theory $\mathcal{CDO}_X$ whose observables
is chiral differential operators on $X$
is holomorphic Chern-Simons theory 
(a.k.a. a curved $\beta\gamma$ system) \cite{GGW}.
It is a 2-dimensional field theory.

An invariant that we have not talked about for field theories
is called the \emph{partition function}.
The partition function of $\mathcal{CDO}_X$ 
recovers an invariant of $X$ known as the Witten genus.
This appears in \cite{chiral,GGW,Costello}.

The Witten genus has close ties to the geometry of elliptic curves.
In fact, there is a cohomology theory built from elliptic curves,
called $tmf$ \cite[Ch. 8]{tmfbook},
that encodes the Witten genus as an orientation. 
An important aspect of $tmf$ is that is has 
\emph{chromatic height 2}. 
That is,
it is part of a hierarchy of more and more complex cohomology theories
of increasing chromatic height.

Similarly, the associative algebra analogue of chiral differential operators 
is differential operators. 
Differential operators on $X$ are related to the observables
of a 1-dimensional theory $\mathcal{D}_X$
called 1-dimensional Chern-Simons theory.
The partition function of $\mathcal{D}_X$ recovers the $\hat{A}$-genus of $X$,
\cite{GwilliamGrady}.
The $\hat{A}$-genus is encoded in ordinary cohomology. 
Ordinary cohomology has chromatic height 1.

\begin{center}
\begin{tabularx}{0.8\textwidth} { 
  | >{\raggedright\arraybackslash}X 
  | >{\centering\arraybackslash}X 
    | >{\centering\arraybackslash}X 
      | >{\centering\arraybackslash}X 
  | >{\raggedleft\arraybackslash}X | }
 \hline
 \begin{center} dimension \end{center} & \begin{center}Field Theory \end{center} & \begin{center}Partition Function \end{center} & \begin{center}Cohomology Theory \end{center} & \begin{center}Chromatic Height\end{center}\\
 \hline
\begin{center}1\end{center}   & \begin{center} 1d Chern-Simons theory \end{center}   & \begin{center} related to the $\hat{A}$ genus \end{center}  & \begin{center} ordinary cohomology $H^\bullet(-)$ \end{center} & \begin{center}1\end{center} \\
 \hline
\begin{center}2\end{center}   & \begin{center} holomorphic Chern-Simons theory \end{center} & \begin{center} related to the Witten genus\end{center} & \begin{center}elliptic cohomology $tmf$ \end{center}  & \begin{center}2 \end{center}\\
\hline
\end{tabularx}
\end{center}

In summary, we have examples of a relationship between 
the dimension of a field theory 
and the chromatic height of the cohomology theory in which the partition function is encoded.

A conjecture of Stolz and Teichner proposes a deeper relationship 
between chromatic height and dimension of field theories;
see \cite{ST1,ST2}.

For recent progress on this program see \cite{WG} and other works of Berwick-Evans.

\section{Duality}

We understand field theories by studying their 
observables.
This translation is great for a few things:

\begin{itemize}
\item a precise definition of topological field theory
\item take advantage of the structure of factorization homology
\item quantization algebraically becomes deformationo
\end{itemize}

I want to end by talking about another benefit of this viewpoint.
Basically in all fields of math I think about, 
notions of duality are super exciting. 
For example, Poincar\'e duality in manifold theory.
In your problem sessions, 
you saw a version of this called non-abelian Poincar\'e duality.

\begin{quest}
Are there notions of dualities in field theory and factorization algebras?
\end{quest}

We are going to start on the algebra side.
The cool type of duality for algebras is called \emph{Koszul duality}.

\begin{defn}
Let $A$ be an associative algebra.
The \emph{Koszul dual} of $A$ is the linear dual
\[\bb{D}(A)=(\bb{1}\otimes_A\bb{1})^\vee.\]
Here $\bb{1}$ denotes the trivial (left or right) $A$-module.
\end{defn}

We can recover this construction using factorization homology.
Indeed, 
there is an equivalence
\[\int_{\bb{D^1}}A=\bb{1}\otimes_A\bb{1}\]
for any $\bb{E}_1$-algebra $A$.

This motivates a definition of Koszul duality for $\bb{E}_n$-algebras.

\begin{defn}[Ayala-Francis]
Let $A$ be an $\bb{E}_n$-algebra.
The \emph{Koszul dual} of $A$ is 
\[\bb{D}(A)=\left(\int_{\bb{D}^n}A\right)^\vee.\]
\end{defn}

This is explained in \cite[Thm. 3.3.2]{Zero}.

\begin{rmk}
There is a different original definition of the Koszul dual 
of an $\bb{E}_n$-algebra, due to Ginzburg-Kapranov and Lurie 
if different contexts. 
Really Ayala-Francis' result is that the definition I 
gave above agrees with these previous definitions.
\end{rmk}

Recently, Ching and Salvatore \cite{ChingSalvatore}
proved a long standing conjecture
regarding the Koszul dual of the \emph{operad} $\bb{E}_n$.

\begin{thm}[Ching-Salvatore,  2020]
The $\bb{E}_n$-operad is Koszul dual to itself.
\end{thm}

This was previously known at the level of chain complexes;
see \cite{GJ}.
Ching and Salvatore's result aslo recovers a previously known result on the level of algebras.

\begin{cor}
The Koszul dual of an $\bb{E}_n$-algebra is an $\bb{E}_n$-algebra.
\end{cor}

This can be found several places, 
including \cite[\S5.2.5]{HA}.

\begin{ex}
Let $\mathfrak{g}$ be a Lie algebra.
The Koszul dual of the enveloping algebra $U(\mathfrak{g})$
is the Lie algebra cochains,
\[\bb{D}(U_n\mathfrak{g})=C^\bullet_\mathrm{Lie}(\mathfrak{g}).\]
This is also true for the $\bb{E}_n$-enveloping algebra $U_n(\mathfrak{g})$.
The Koszul dual of $U_n(\mathfrak{g})$ is $C^\bullet_\mathrm{Lie}(\mathfrak{g})$,
viewed as an $\bb{E}_n$-algebra,
see \cite[Cor. 4.2.1]{PKD}.
\end{ex}

Thus, you can ask if the factorization homology 
of dual algebras are related.
The following is the main theorem of \cite{PKD}.

\begin{thm}[Poincar\'e/Koszul duality, Ayala-Francis]
Under conditions,
there is an equivalence
\[\int_MA\simeq\left(\int_{M^+}\bb{D}(A)\right)^\vee.\]
\end{thm}

The changing of $M$ to $M^+$ reflects,
for example, the difference in Poincar\'e duality for manifolds with 
boundary.

Turning back to the field theory side,
note that
the Koszul dual of observables on $\bb{R}^n$ 
has the right algebraic structure to be observables of
a different field theory on $\bb{R}^n$.
That is, 
let $\mathsf{Obs}_\bb{X}$ denote the observables of a field theory $\bb{X}$.
We would like to know if there is another field theory $\bb{Y}$ so that 
the Koszul dual of observables on $\bb{X}$ is observables on $\bb{Y}$;
in symbols
\[\bb{D}(\mathsf{Obs}_{\bb{X}})\simeq \mathsf{Obs}_{\bb{Y}}.\]

\begin{conj}[Costello-Paquette, Costello-Li]
The holomorphically twisted version of AdS/CFT duality
gives an example of Koszul duality for observables.
\end{conj}
See \cite{CP} and the references therein.

To make this conjecture precise,
one would need a notion of Koszul duality for factorization algebras.
That is an open problem.
Checking this conjecture in examples is an active 
area of research by lots of people.

\section{Exercises}

Let $Z$ be a fully extended $n$-dimensional TQFT
valued in $\mathcal{C}$ with 
$\Omega^n\mathcal{C}=\bb{C}$,
$\Omega^{n-1}\mathcal{C}=\mathsf{Ch}_\bb{C}$,
and
$\Omega^{n-2}\mathcal{C}=\mathsf{LinCat}_\bb{C}$.

Recall that 
$Z(S^{n-1})$ is an $\bb{E}_n$-algebrain $\mathsf{Ch}$
$Z(S^{n-2})$ is an $\bb{E}_{n-1}$-algebra in $\mathsf{LinCat}$.

\begin{quest}\label{q-20.1}
Show the analogous statement for $Z(S^{n-k})$.
\end{quest}

\begin{quest}\label{q-20.2}
Let $K$ be a $(n-1)$-manifold with boundary.
Show that $K$ determines an object in $Z(\partial K)$.
\end{quest}

\begin{quest}\label{q-20.3}
Can you describe the unit of $Z(S^{n-2})$?
\end{quest}

\begin{quest}\label{q-20.4}
Let $A$ be an $\bb{E}_n$-algebra. 
Recall from yesterday's exercises that you constructed
an $(\infty,n)$-category from $A$. 
Use this construction to write down sensible values 
for an $n$-dimensional field theory $Z$
with $Z(S^{n-1})=A$.
In particular, write down $Z(S^{n-k})$ for all $k$,
and $Z(\mathrm{pt})$.
\end{quest}

\begin{quest}\label{q-20.5}
Show that a vertex algebra $(V,Y,T)$ with 
\[Y(u,z)\in\mathsf{End}V[[z]]\]
is equivalent to one formed by a commutative ring with derivation.
\end{quest}

\newpage
\part{Solutions and Hints for Exercises}

\subsubsection{Exercise \ref{q-5.new0}}
See \cite[Part 1, \S 3]{Curry}.

\subsubsection{Exercise \ref{q-5.2}}
The open cover consisting of $U_1=(-1,\infty)\times\bb{R}$
and $U_2=(-\infty,1)\times\bb{R}$ is not a Weiss cover.
The finite set $\{(-6,0),(6,0)\}$ is not fully contained in any open set of the cover.

\subsubsection{Exercise \ref{q-5.3}}
Take the open cover consisting of just $M$ itself.
More interestingly, you could take the open cover
consisting of all disjoint unions of embedded open disks in $M$. 

\subsubsection{Exercise \ref{q-5.5}}
An open cover of $\mathsf{Ran}(M)$ 
is a collection of open sets $\{U_i\}$
so every finite set of points in $M$ is contained in some $U_i$.

\subsubsection{Exercise \ref{q-1comm}.}
Let $A$ be a commutative algebra (say in chain complexes). 
The constant precosheaf defines a functor 
\[F_1\colon\mathsf{Open}(M)\rta\mathsf{Ch}\]
sending every open set to $A$.
We can view this as a precosheaf valued in algebras in chain complexes,
\[F_2\colon\mathsf{Open}(M)\rta\mathsf{Comm}(\mathsf{Ch}).\]
We can cosheafify $F_2$ using the Weiss topology to get a Weiss cosheaf
\[F_3\colon\mathsf{Open}(M)\rta\mathsf{Comm}(\mathsf{Ch}).\]
Since coproducts and tensor products are the same in algebras in chain complexes,
we have 
\[F_3(U\sqcup V)\simeq F_3(U)\otimes F_3(V).\]
Forgetting back down,
we have a factorization algebra
\[F_4\colon\mathsf{Open}(M)\rta\mathsf{Ch}.\]
The functor $F_4$ is still a Weiss cosheaf since the forgetful functor 
\[\mathsf{Comm}(\mathsf{Ch})\rta\mathsf{Ch}\]
preserves reflexive coequalizers.

\subsubsection{Exercise \ref{q-5.7}}
Let $A$ be an associative algebra.
This should be a factorization algebra on $\bb{R}$. 
It assigns the associative algebra $A$ to every open interval.
The factorization algebra structure comes from the multiplication on $A$.

\subsubsection{Exercise \ref{q-5.8}}
When the factorization algebra $\mathcal{F}$
sends every inclusion $(a,b)\subset (c,d)$ to
an equivalence.

\subsubsection{Exercise \ref{q-5.9}}
Your answer should involve right and left modules.

\subsubsection{Exercise \ref{q-5.new}}
The space of fields is a vector space,
saw $W=\bb{R}^n$. 
The critical locus of $S$ is then the critical locus of a function 
\[S\colon\bb{R}^n\rta \bb{R}.\]

\subsubsection{Exercise \ref{q-5.10}}
The critical locus is the set of points with $dS(p)=0$. 
The zero section of $T^*Y$ is the set of pairs $(p,0)$.
The graph of $dS$ is the set of pairs $(p, dS(p))$.

\subsubsection{Exercise \ref{q-5.11}}
If $S=0$, then $dS=0$.

\subsubsection{Exercise \ref{q-5.13}}
See \cite[\S 1]{Nat}.

\subsubsection{Exercise \ref{q-5.14}}
See \cite[\S 1]{Nat}.

\subsubsection{Exercise \ref{q-10.1}}
An $\bb{E}_n$-algebra is an $\bb{E}_m$-algebra for $m<n$.
If you have already done Exercise \ref{q-10.8},
you can see this geometrically by May's recognition principle.
The space $\Omega^nX$ is the same as
$\Omega^mY$ where $Y=\Omega^{n-m}X$.

\subsubsection{Exercise \ref{q-10.2}}
Draw pictures of disjoint unions of disks including into other disks.

\subsubsection{Exercise \ref{q-10.3}}
See \cite[Ex. 1.1.7]{LurieEk}.

\subsubsection{Exercise \ref{q-10.4}}
See \cite[Cor. 1.1.9]{LurieEk}.
Alternatively, one can see this on the level of operads
as is done in \cite[Thm. 4.8]{May}.

\subsubsection{Exercise \ref{q-10.6}}
See \cite[Ex. 1.2.4]{LurieEk}.
Alternatively, one can see this on the level of operads
as is done in \cite[Thm. 4.8]{May}.

\subsubsection{Exercise \ref{q-10.7}}
See \cite[Ex. 1.2.4]{LurieEk}.

\subsubsection{Exercise \ref{q-10.5}}
See \cite[Cor. 1.1.9]{LurieEk}.
Alternatively, one can see this on the level of operads
as is done in \cite[Thm. 4.8]{May}.

\subsubsection{Exercise \ref{q-10.8}}
Think of $\Omega^nX$ as 
compactly supported maps $\bb{R}^n\rta X$ 
and using inclusion of disjoint unions of disks into other disks.
See \cite[Thm. 5.1]{May}.

\subsubsection{Exercise \ref{q-10.9}}
Start by checking that disjoint union of disks is sent to a tensor product.

\subsubsection{Exercise \ref{q-10.10}}
See \cite[Ch. 3 \S 4]{CG1}.

\subsubsection{Exercise \ref{q-10.11}}
See \cite[Ch. 3 \S 4]{CG1}..

\subsubsection{Exercise \ref{q-12.1}}
Use excision and the description of a factorization algebra on $[0,1]$ from a previous exercise.

\subsubsection{Exercise \ref{q-12.2}}
Use excision.

\subsubsection{Exercise \ref{q-3excision}}
Here $F(V\times\bb{R})$ inherits an 
$\cal{E}_1$-algebra structure from the copy of $\bb{R}^1$,

\[\msf{Emb}^\mrm{fr}(\coprod_I\bb{R},\bb{R})\otimes F(V\times\bb{R})^{\otimes I}\simeq \msf{Emb}^\mrm{fr}(\coprod_I\bb{R},\bb{R})\otimes F(V\times(\coprod_I\bb{R}))\rta F(V\times\bb{R})\]

The tensor product

\[F(U)\bigotimes_{F(V\times\bb{R})}F(U')\rta F(W)\]
is then the tensor product in modules over the $\cal{E}_1$-algebra $F(V\times\bb{R})$. 

\subsubsection{Exercise \ref{q-12.4}}
Show that $C_\bullet $ takes $\bb{E}_n$-algebras in spaces 
to $\bb{E}_n$-algebras in chain complexes.

\subsubsection{Exercise \ref{q-12.5}}
The proof is given in \cite[Lem. 4.5]{AFtop}.

\subsubsection{Exercise \ref{q-12.6}}
Use the fact that every homology theory for manifolds
looks like factorization homology for some $\bb{E}_n$-algebra.

\subsubsection{Exercise \ref{q-12.7}}
See \cite[Ch. 3 \S 4]{CG1}.

\subsubsection{Exercise \ref{q-12.8}}
See \cite[Ch. 3 \S 4]{CG1}.

\subsubsection{Exercise \ref{q-3Ug}}
Use both the fact that $C_\bullet^\mathrm{Lie}$ commutes with 
factorization homology
and nonabelian Poincar\'e duality.

\subsubsection{Exercise \ref{q-16.1}}
In general, we should take $Z(S^{n-1})$.
Here the dimension is $n=1$ so we are looking at $Z(S^0)$.
Since $Z_V$ sends a positively oriented point to $V$ 
and a negatively oriented point to $V^\vee$, 
we have $Z(S^0)=\mathsf{End}(V)$.

\subsubsection{Exercise \ref{q-16.2}}
Consider the bordism $N\times[0,1]$.
Remove an open disk from $N\times\{1/2\}$.
This creates a bordism from $N\sqcup S^{n-1}$ to $N$.

\subsubsection{Exercise \ref{q-16.3}}
Let $\rho\colon G\rta \mathsf{Aut}(W)$ be the representation.
Given a principal $G$-bundle $P\rta M$, 
build the associated bundle $P\times _GW$.
Take the holonomy as when constructing Wilson loop operators.
The fiber of $P\times_GW$ is $W$
so the holonomy gives a map $W\rta W$.
Take the trace of this map as a linear map.

\subsubsection{Exercise \ref{q-16.5}}
One category is $\mathsf{Mod}_A$ of modules.
Another is $BA$ with a single object and morphisms $A$.

\subsubsection{Exercise \ref{q-16.6}}
Oriented closed $2$-manifolds, all of which look like 
disjoint unions of spheres and multi-hole torii.
The empty set is also an object.

\subsubsection{Exercise \ref{q-16.7}}
This is the classification or oriented surfaces.

\subsubsection{Exercise \ref{q-16.8}}
We can think of building the two-holed torus
as a connect sum as analogous to splitting the 
two-holed torus into a composition of morphisms in $\mathsf{Cob}(2)$.

\subsubsection{Exercise \ref{q-16.9}}
Read the next part for the answer!

\subsubsection{Exercise \ref{q-20.1}}
Use the pair of pants bordisms of dimension $n-k+1$. 
That is, the bordism from multiple disjoint copies of $S^{n-k}$
to a single copy of $S^{n-k}$.

\subsubsection{Exercise \ref{q-20.2}}
View $K$ as a morphism in the bordism category from the empty set to $\del K$.
The functor $Z$ takes this morphism to a morphism 
from the unit $\bb{C}$-linear category $\mathsf{Ch}_\bb{C}$ to $Z(\del K)$. 
The object $\bb{C}$ is sent to an object in $Z(\del K)$.
This is the object $K$ correspond to.

\subsubsection{Exercise \ref{q-20.3}}
Consider the object coming from the bordism $\bb{D}^{n-1}$.

\subsubsection{Exercise \ref{q-20.4}}
Consider the $(\infty,k)$-categories $k\mathsf{Mod}_A$ 
of iterated modules.

\subsubsection{Exercise \ref{q-20.5}}
See \cite[Ch. 1 \S 1.4]{BenZviFrenkel}.

\newpage
\bibliography{bib}
\bibliographystyle{alpha}

\end{document}